\documentclass[review]{siamart190516}

\usepackage{braket,amsfonts,soul}

\usepackage{array}

\usepackage[caption=false]{subfig}
\usepackage{pgfplots}

\newsiamthm{claim}{Claim}
\newsiamremark{remark}{Remark}
\newsiamremark{hypothesis}{Hypothesis}
\newsiamremark{test}{Test}
\crefname{hypothesis}{Hypothesis}{Hypotheses}
\newsiamremark{exercise}{Exercise}
\newsiamremark{example}{Example}
\usepackage{algorithmic}

\usepackage{graphicx,epstopdf}

\Crefname{ALC@unique}{Line}{Lines}

\usepackage{amsopn}

\usepackage{xspace}
\usepackage{bold-extra}
\usepackage[most]{tcolorbox}

\newcounter{example}\stepcounter{example}
\colorlet{texcscolor}{blue!50!black}
\colorlet{texemcolor}{red!70!black}
\colorlet{texpreamble}{red!70!black}
\colorlet{codebackground}{black!25!white!25}

\newcommand{\di}{\partial}
\newcommand{\exsprucebudworm}{{\it Spruce Budworm.} }
\newcommand{\exflexiblebeam}{{\it Flexible Beam.} }
\newcommand{\expendulum}{{\it Pendulum.} }
\newcommand{\exferro}{{\it Ferromagnetism.} }

\newcommand{\R}{ \mathbb R} 

\newcommand{\achg}[1]{{\color{blue} #1}}

\newcommand{\kchg}[1]{{\color{red}  #1}}

\begin{tcbverbatimwrite}{tmp_\jobname_header.tex}
\nolinenumbers
\title{Hysteresis and Stabillity \thanks{Submitted to the editors \today
\funding{The financial support  of  NSERC for the research discussed in this article is gratefully acknowledged.}}}

\author{Amenda Chow\thanks{Department of Mathematics and Statistics, York University, Canada (\email{amchow@yorku.ca}).}
\and Kirsten A. Morris\thanks{Department of Applied Mathematics, University of Waterloo, Canada (\email{kmorris@uwaterloo.ca}).}
\and Gina Faraj Rabbah\thanks{Department of Mathematics, University of Toronto, Canada (\email{gnicole98@gmail.com})}}

\end{tcbverbatimwrite}
\nolinenumbers
\title{Hysteresis and Stabillity \thanks{Submitted to the editors \today
\funding{The financial support  of  NSERC for the research discussed in this article is gratefully acknowledged.}}}

\author{Amenda Chow\thanks{Department of Mathematics and Statistics, York University, Canada (\email{amchow@yorku.ca}).}
\and Kirsten A. Morris\thanks{Department of Applied Mathematics, University of Waterloo, Canada (\email{kmorris@uwaterloo.ca}).}
\and Gina Faraj Rabbah\thanks{Department of Mathematics, University of Toronto, Canada (\email{gnicole98@gmail.com})}}

\ifpdf
\hypersetup{ pdftitle={SIAM Review hysteresis article} }
\fi

\begin{document}
\maketitle
\nolinenumbers
\begin{tcbverbatimwrite}{tmp_\jobname_abstract.tex}
\begin{abstract}
 Hysteresis can be  defined from a dynamical systems perspective with respect to  equilibrium points. Consequently, hysteresis naturally lends itself as a topic to illustrate  and extend concepts in a dynamical systems course. A number of examples exhibiting hysteresis, most motivated by applications, are presented. Although the examples can be used to construct student exercises,  specific questions are listed in an appendix.  A brief extension on hysteresis in partial differential equations is also included.  
 \end{abstract}


\end{tcbverbatimwrite}
\begin{abstract}
 Hysteresis can be  defined from a dynamical systems perspective with respect to  equilibrium points. Consequently, hysteresis naturally lends itself as a topic to illustrate  and extend concepts in a dynamical systems course. A number of examples exhibiting hysteresis, most motivated by applications, are presented. Although the examples can be used to construct student exercises,  specific questions are listed in an appendix.  A brief extension on hysteresis in partial differential equations is also included.
 \end{abstract}




\section{Introduction}\label{secintro}
Hysteresis appears in many contexts; such as biology \cite{Angeli2004,Collings1990,Mayergoyz2019,Noori2014}, magnetism \cite{Chen2019,Chow2014_ACC, Chow2015,Bernstein2005}, ferroelectric materials \cite{Xu2016,Park2019,Wang2009}, climate change \cite{Berdugo2021}, circuits \cite{Elwakil2013,Sharifi2015} and economics \cite{Grinfield2009}.  Many researchers have highlighted the relationship between multiple stable equilibria and hysteresis in various applications, including \cite{Afshar2016,Alimov1998,Angeli2004,Sontag2003,Campbell2021,Chow2014_ACC,Collings1990,Huang2018,Ikhouane2013,Morris2011,Padthe2008,Sharifi2016,Wallis1994,Wang2017}.
The looping behaviour characteristic of hysteresis can be explained through analysis of equilibria and stability; see for example, the review article \cite{Morris2011}.

  A dynamical systems viewpoint provides a fundamental perspective on the cause of hysteresis. Consequently, the teaching of dynamical systems concepts can be reinforced through  hysteresis. In the next section,  hysteresis is defined formally and its connection to stability is described. This is illustrated by examples of ordinary differential equations  (ODEs) that model  various applications. Simulations are presented to  complement the analysis.    A brief discussion of possible extensions is provided, as well as a number of problems suitable for students. 

\section{Stability of Equilibria}\label{secdefining}
Consider a general system of ordinary differential equations with state $\mathbf x (t) \in \R^n $ and  a real-value continuous forcing function $\mathbf  u (\cdot);$
\begin{equation}
\label{FONODE}
\dot {\mathbf x}(t)= \mathbf f(\mathbf x(t),\mathbf u(t)), \quad   t \geq 0, \quad
\mathbf x(0) =\mathbf  x_0,
\end{equation}
where $\mathbf f$ is a continuously differentiable function of  $\mathbf  x$,  and is continuous with respect to $\mathbf  u.$ We call $\mathbf  u$ the input.

The following definitions are standard, see for example \cite[page~112]{Khalil2002},
\begin{definition}
If $f(\bar{\mathbf  x})=0$, then $\bar {\mathbf x}$ is an {\em equilibrium point} of \eqref{FONODE}.  

An equilibrium point of \eqref{FONODE} is {\em stable} if for each $\epsilon>0$ there is a $\delta>0$ such that $||\mathbf x(0)||<\delta$ implies $||\bar {\mathbf x}-{\mathbf  x}(t)|| <\epsilon$ for all $t\geq0. $  If an equilibrium point is not stable, it is {\em unstable}.  
If there is a $\delta>0$ such that $||\mathbf x(0)||<\delta$ implies $\displaystyle\lim_{t \to \infty } ||\bar {\mathbf x}- {\mathbf x}(t)|| =0 $ the equilibrium point is {\em asymptotically stable}.
\end{definition}
For different constant inputs, the equilibria of the  differential equation in \eqref{FONODE} will, in general, depend on the input. These equilibria are fundamental to the existence and analysis of hysteresis.  The simplest method to analyze stability of equilibia is through Lyapunov's Indirect Method. Denote the Jacobian of $\mathbf  f$ with respect to $\mathbf x$ by $D \mathbf f$ with $D \mathbf f$ evaluated at an equilibrium point of \eqref{FONODE}.  If the real part of the eigenvalues of $D\mathbf f$ are all negative, then  the equilibrium point is asymptotically stable; if the real part of at least one eigenvalue of  $D\mathbf f$ is positive,  then the equilibrium point is unstable; otherwise, this approach yields no conclusion on stability. For additional details of Lyapunov's Indirect Method, see \cite[section~4.3]{Khalil2002}.

\begin{example}\label{exmotivation}
Consider the following ODE
\begin{align}\label{exnonlinearODEPerko}
\dot x(t)& =x(t)(1-(x(t))^2) + u(t).
 \end{align}

Equation~\eqref{exnonlinearODEPerko} is a  differential equation analyzed in many texts; see for example,  \cite[Page 336]{Perko2002}, \cite{Bernstein2009}, \cite[Page 329]{Berglund2000}. In \cite[Page 336]{Perko2002},  it is used to illustrate bifurcation at equilibria.

The equilibrium points of  \eqref{exnonlinearODEPerko} for constant scalar input $u(t)=U$  where $U \in \R,$ are  the roots of 
\[
x - x^3 + U=0.
\]
 Figure \ref{PerkoBifurcation} shows these roots  as a function of $U$.    For $U \in (-0.385,0.385)$, there are three equilibria with one unstable equilibrium point between two stable equilibria and the stable equilibria are equidistant from the unstable equilibrium; and, outside of this range there is a single, stable, equilibrium point.  For instance, when  $U=1$,  equation~\eqref{exnonlinearODEPerko} has one equilibrium point. For $U=0$, equation~\eqref{exnonlinearODEPerko} has three equilibrium points: $1$, $0$ and $-1$. By considering the derivative of $f$ with respect to $x$, the equilibrium points $1$ and $-1$ are shown to be  stable, while $0$ is unstable.

Regarding $u$ as the input and $x$ as the solution to \eqref{exnonlinearODEPerko}, Figure~\ref{fig:Perko5} shows the graph of $x$ as a function of $u$ where $u(t)=\sin(\omega t)$  for different frequencies $\omega.$ Loops are seen, which persist as  the frequency $\omega$ goes to zero. Figure \ref{fig:PerkoZoomed} is a magnified version of one of the loops of Figure~\ref{fig:Perko5} and shows jumps in the value of $x$ when the periodic input  takes values around -0.385 and 0.385. This is where the system moves between the stable equilibria 1 and -1.

\begin{figure} [h] 
\centering
\includegraphics[width = .25 \linewidth]{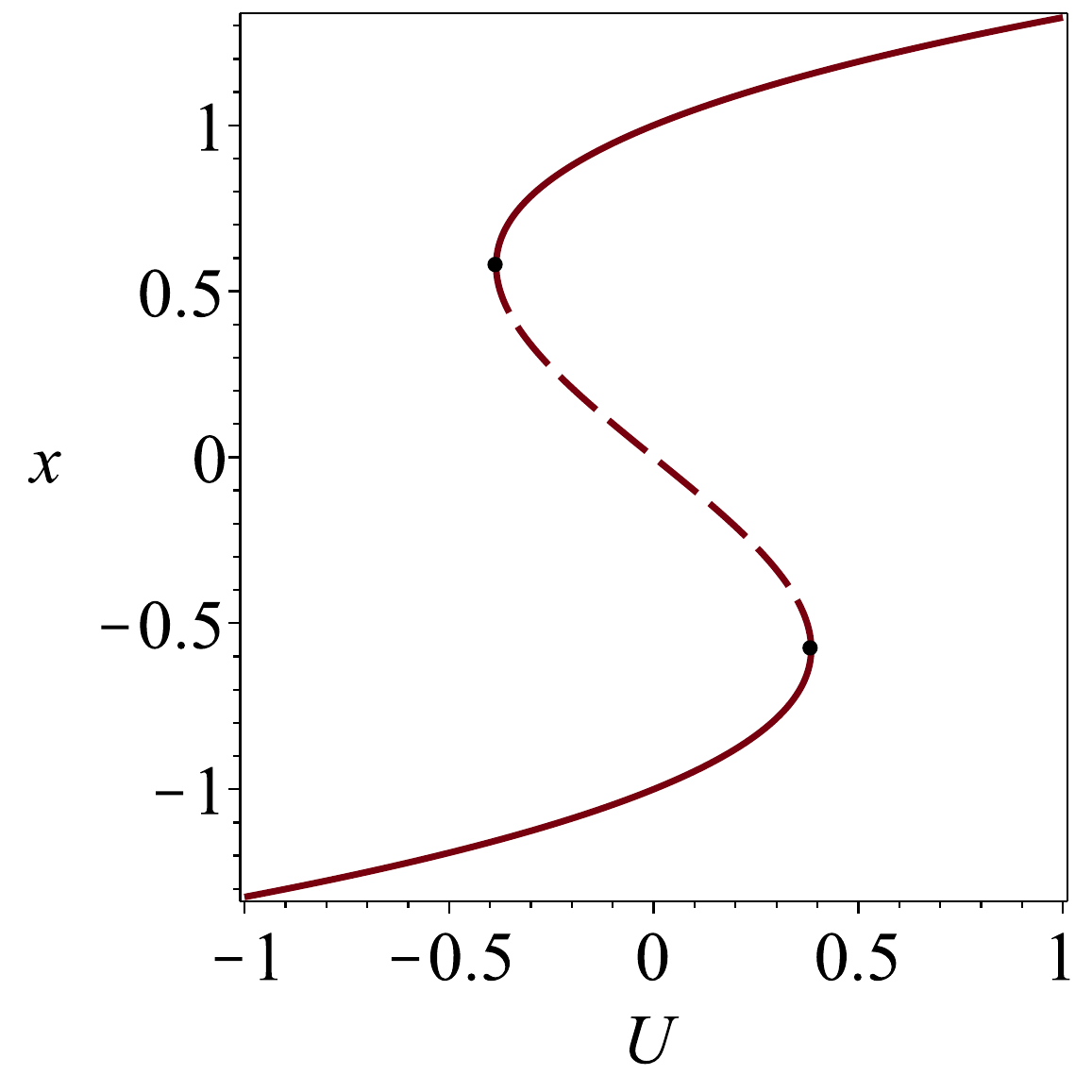}
\caption{\small \label{PerkoBifurcation} This diagram depicts the equilibria of \eqref{exnonlinearODEPerko} as a function of $U$, where the dashed portion indicates  unstable equilibria and the solid portion indicates stable equilibria. This figure was constructed using Maple.}
\end{figure}

\begin{figure}[h!]
\begin{center}
  \subfloat[]{\label{Pa}\includegraphics[width = .24\linewidth, trim=0 175 0 175]{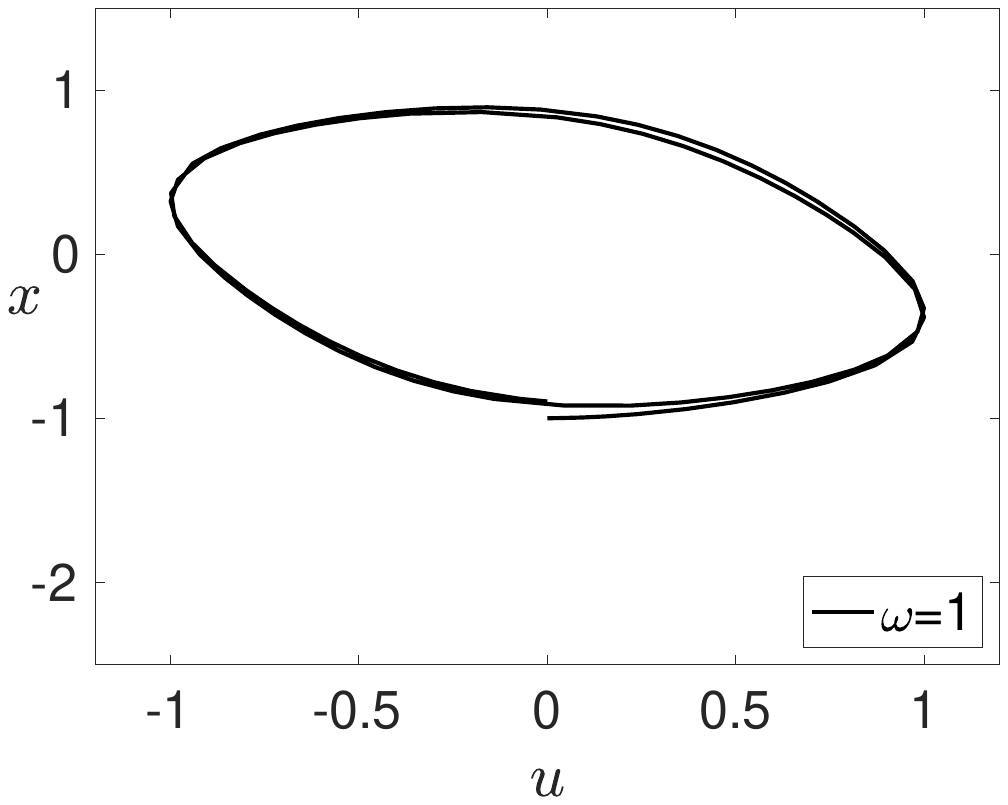}}
 \subfloat[]{\label{Pb}\includegraphics[width = .24 \linewidth, trim=0 175 0 175]{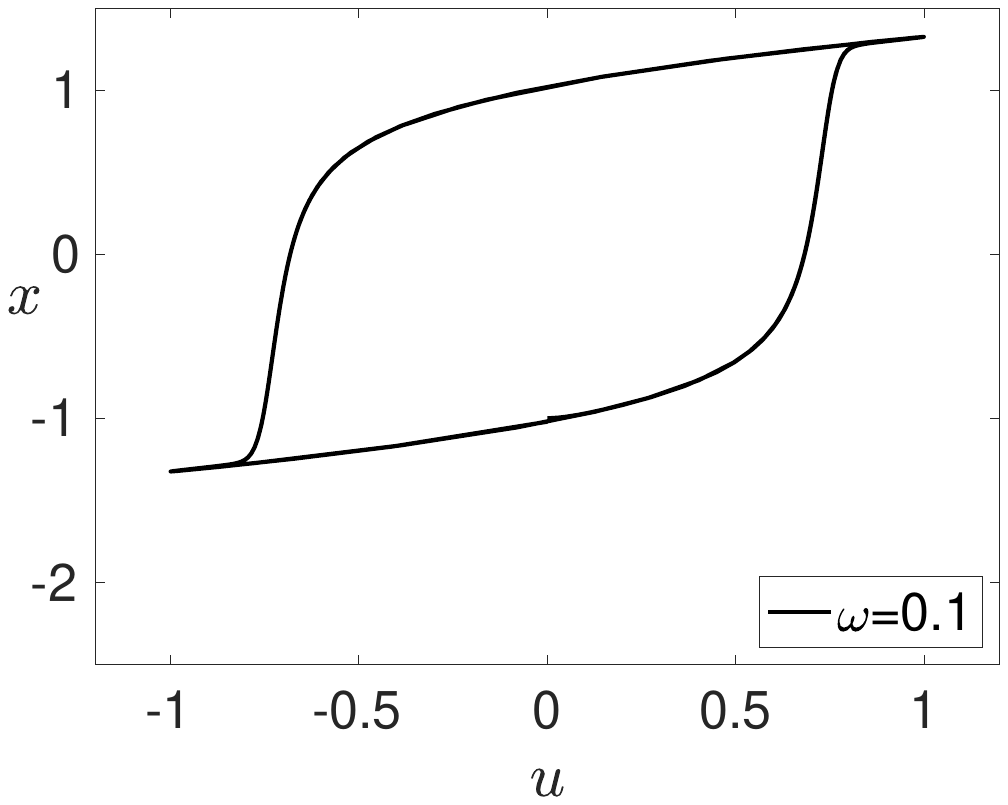}} 
   \subfloat[]{\label{Pc}\includegraphics[width = .24 \linewidth, trim=0 175 0 175]{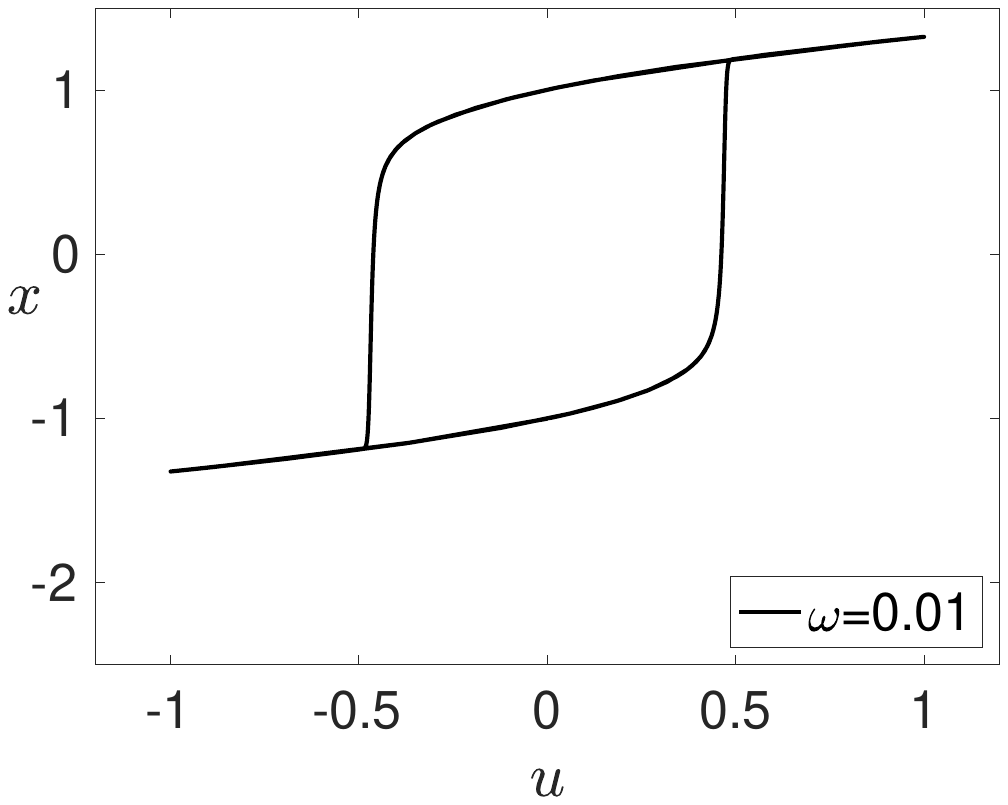}} 
\subfloat[]{\label{Pe}\includegraphics[width = .24 \linewidth, trim=0 175 0 175]{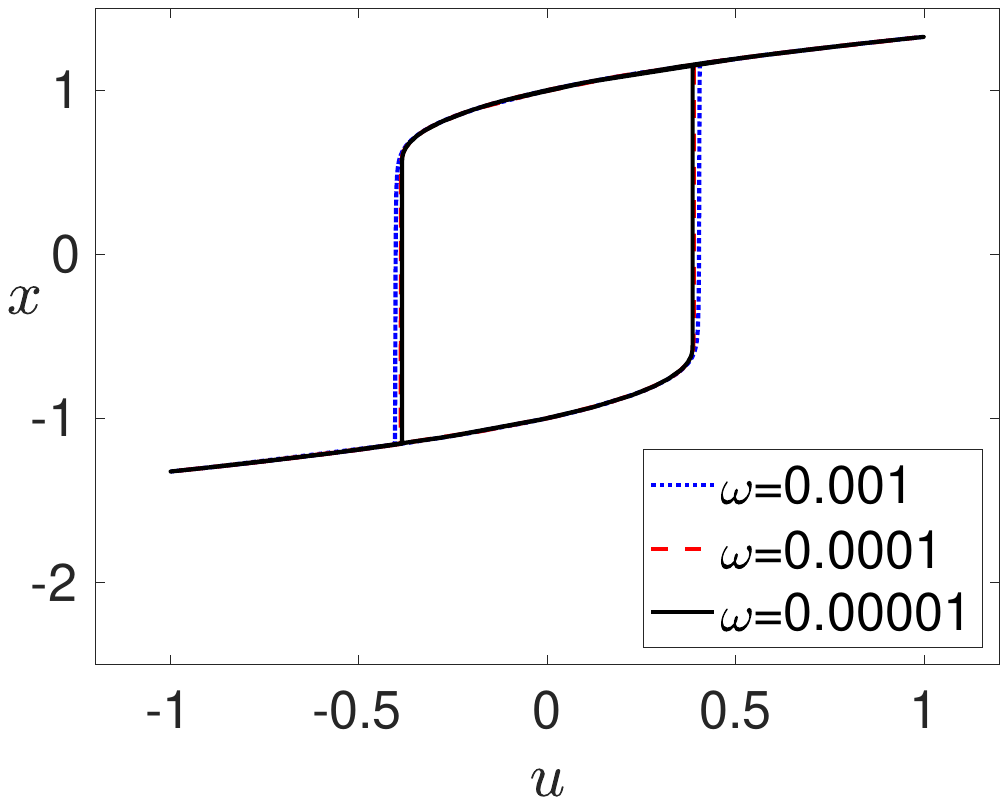}}
\caption{\small Depicted are the graphs of the solution of \eqref{exnonlinearODEPerko}, $x,$ as a function of the input $u(t)=\sin(\omega t)$ for various values of $\omega$. The initial condition chosen is $x(0)=-1$. Equation~\eqref{exnonlinearODEPerko} has two stable equilibria and one unstable equilibrium point and exhibits looping behaviour.}
\label{fig:Perko5}
\end{center}
\end{figure}

\begin{figure}[h!]
\begin{center}
\includegraphics[width = 0.35 \linewidth, trim={0 6cm 0 3cm},clip]{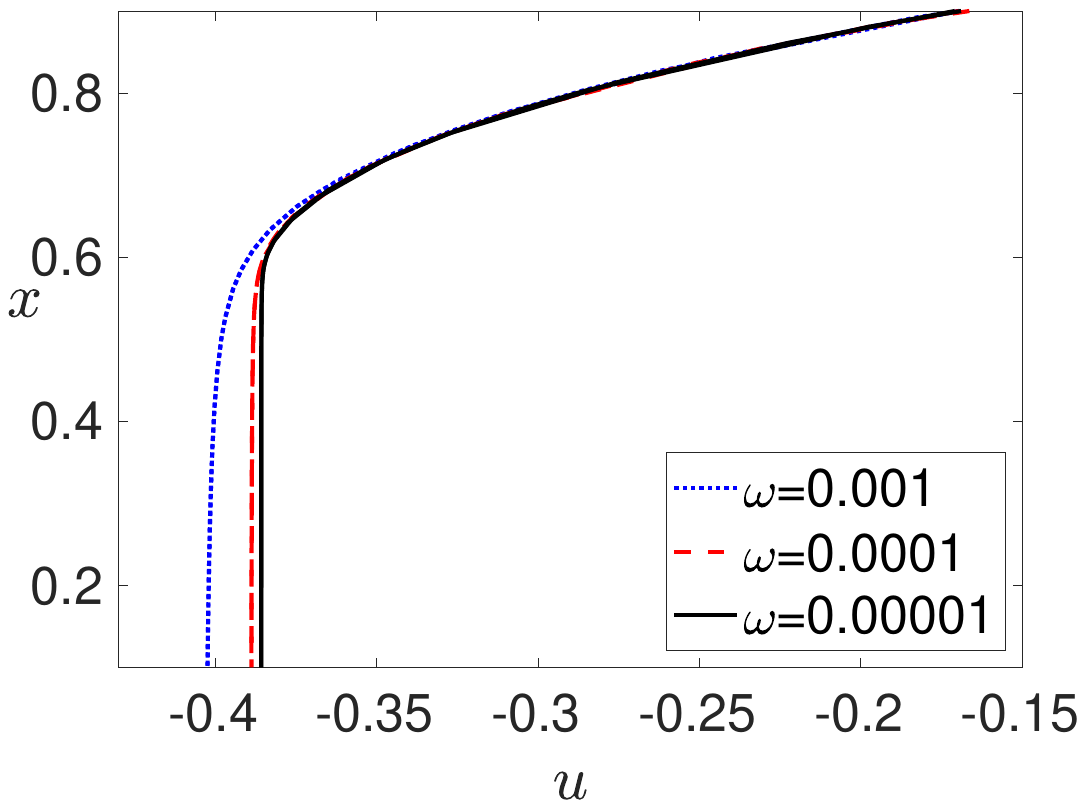}
\includegraphics[width = 0.35 \linewidth, trim={0 6cm 0 3cm},clip]{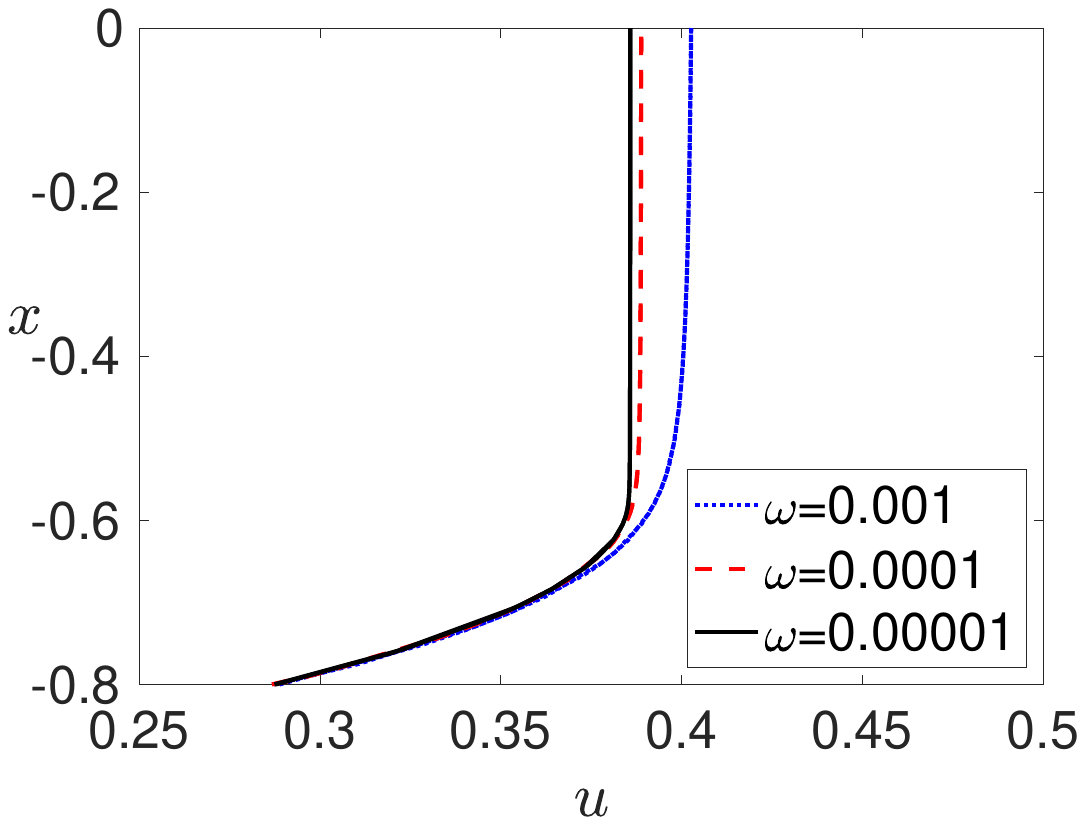}
\caption{\small A magnified image of Figure~\ref{Pe} shows the vertical jumps of the loop converge to $-0.385$ and $0.385$ as $\omega$ approaches 0. \label{fig:PerkoZoomed}}
\end{center}
\end{figure}

\end{example}

\section{Hysteresis}

The looping behaviour and presence of multiple stable equilibria illustrated in Example \ref{exmotivation}  is associated with a phenomenon known as hysteresis.
\begin{definition}\cite[Definition~3]{Morris2011}~\label{defMorris}\\
A system exhibits \textit{hysteretic} if it has\\
 (i) multiple stable equilibria and \\
 (ii) system dynamics that are faster than the time scale at which inputs are varied.
 \end{definition}
   
\begin{figure}[h!]
\begin{center}
\includegraphics[width = 4.1 cm, trim=1cm 6.8cm 1cm 5cm , clip]{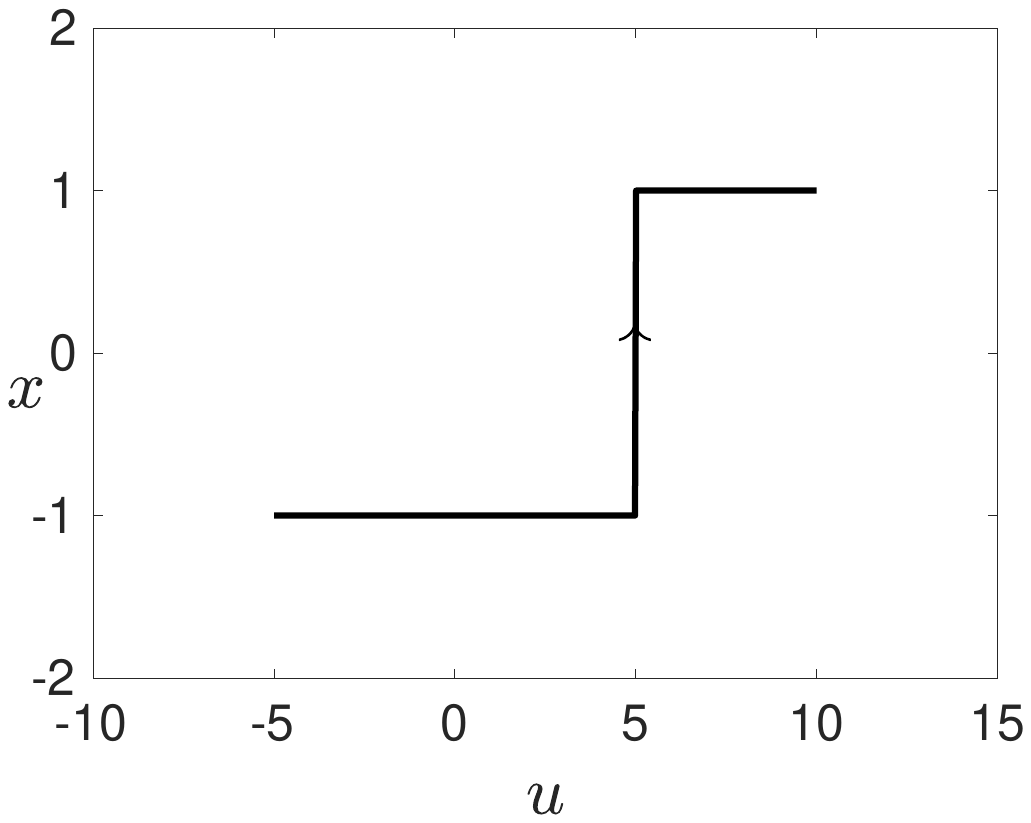} \hspace{0 cm}
\includegraphics[width = 4.1 cm, trim=1cm 6.8cm 1cm 5cm , clip]{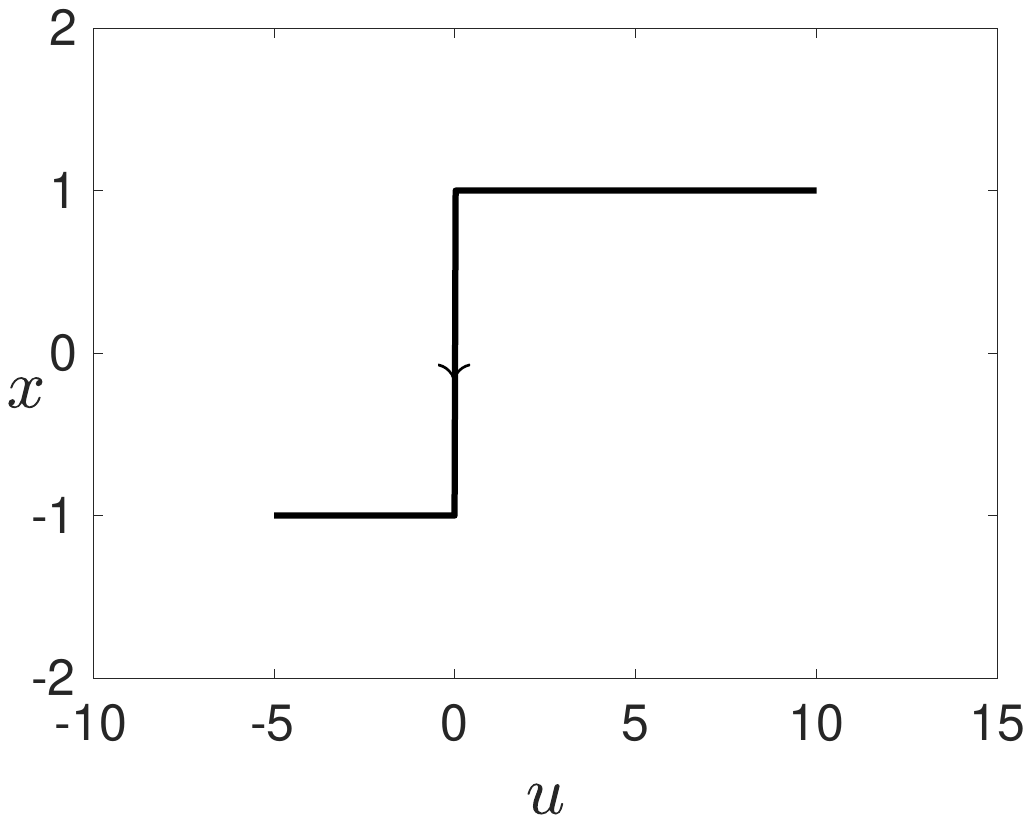} \hspace{0 cm}
\includegraphics[width = 4.1 cm, trim=1cm 6.8cm 1cm 5cm , clip]{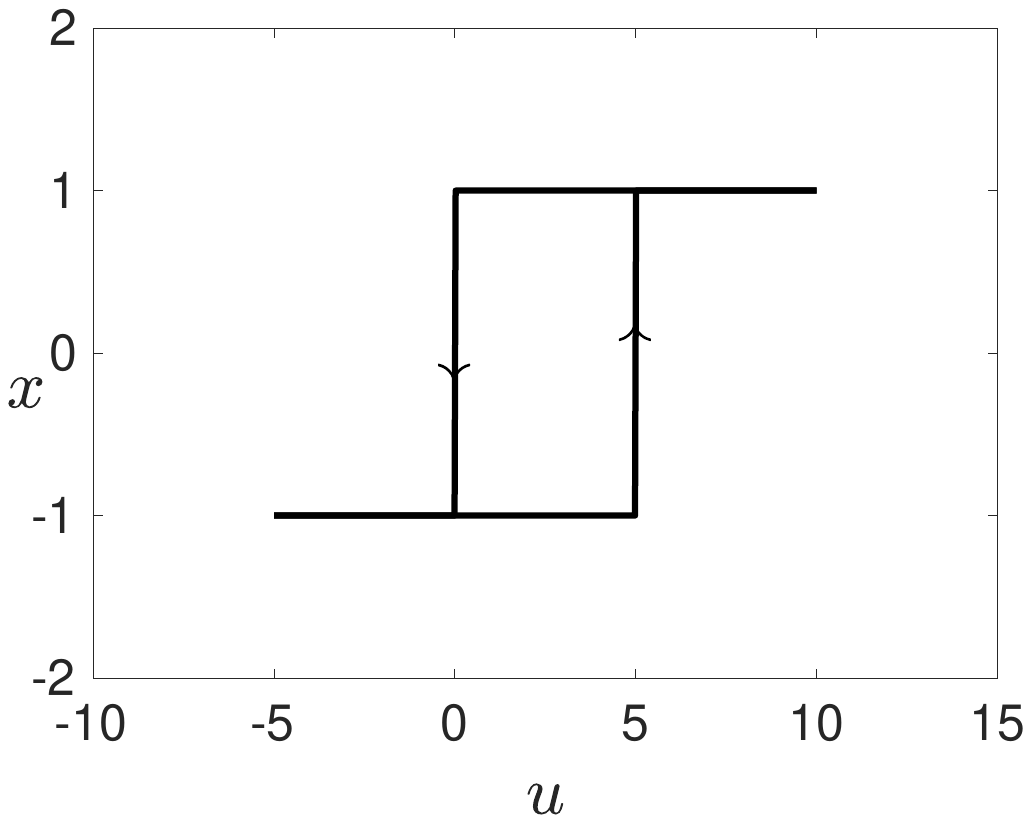}
\caption{\small  Suppose -1 and 1 are two stable equilibria of $\dot x(t)=f(x(t),u(t))$ and $x\in \R$ is its solution. In the left, the system dynamics are depicted as ``jumping" from -1 to 1 as the periodic input $u(t)$ varies. In the middle, as $u(t)$ continues to vary, the system dynamics follows a different path from 1 to -1.  This creates a loop as shown in the plot on the right.}
\label{figrelay}
\end{center}
\end{figure}

As illustration, suppose for \eqref{FONODE} being scalar that  $f(-1,u_1) =0$ and $f(1, u_2)=0 $ with $1$ and $-1$ are   stable equilibria of $f(x(t),u)$ and suppose $u_1 < u_2. $ When $x(t)$ is  near $1$ and  the input $u(t)$ is near $u_2$ the system dynamics tend to stay near that equilibrium.  But if the input varies enough the system state shifts to the basin of attraction of the other stable equilibrium point. If system dynamics are fast compared to the change in the value of the input, this transition between equilibria sometimes appears to occur almost instantaneously.   Due to this transition between equilibria, the path of the system  is different when the input  increases than when it decreases.  This leads to the loops characteristic of hysteresis as shown in Figures \ref{fig:Perko5} and \ref{figrelay}.
   
\begin{figure}[h!]
\begin{center}
\vspace{-1.6cm} \subfloat{\includegraphics[width = 4.3 cm]{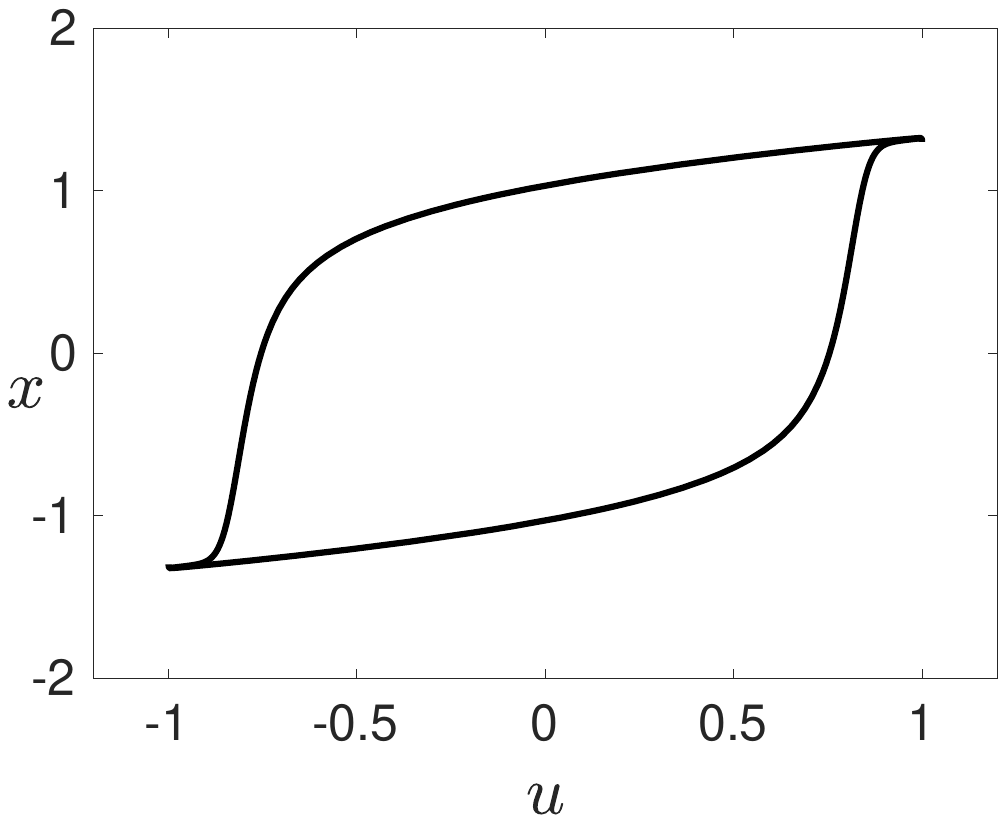}} \hspace{0 cm}
\vspace{-1.6cm} \subfloat{\includegraphics[width = 4.3 cm]{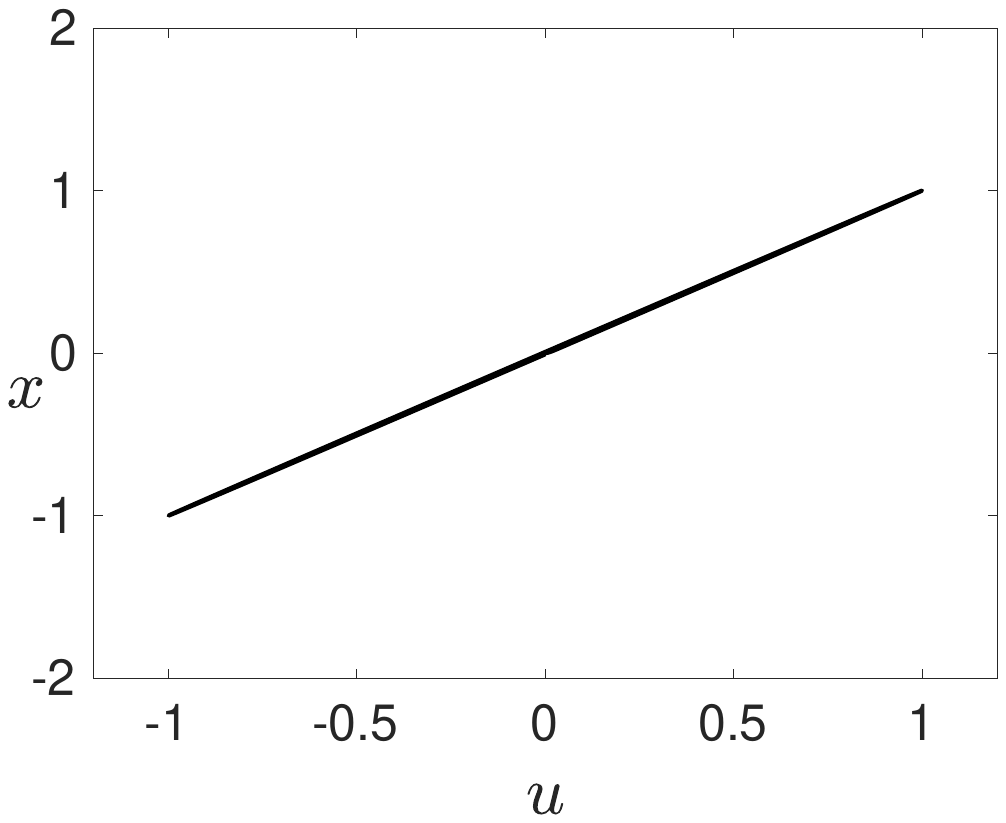}} \hspace{0 cm}\\
\caption{\small   For a single differential equation with solution $x \in \mathbb R$ and periodic input $u \in \mathbb R$, displayed is the graph of $x$ as a function of $u$ with one displaying looping behaviour (nontrivial closed curve, left) and the other not displaying looping (trivial closed curve, right). }
\label{fig:comparehysteresis}
\end{center}
\end{figure}

\begin{figure}[h!]
\begin{center}
\vspace{-1.6cm} \subfloat{\includegraphics[width = 4.3 cm, trim=1cm 6.8cm 1cm 0cm , clip]{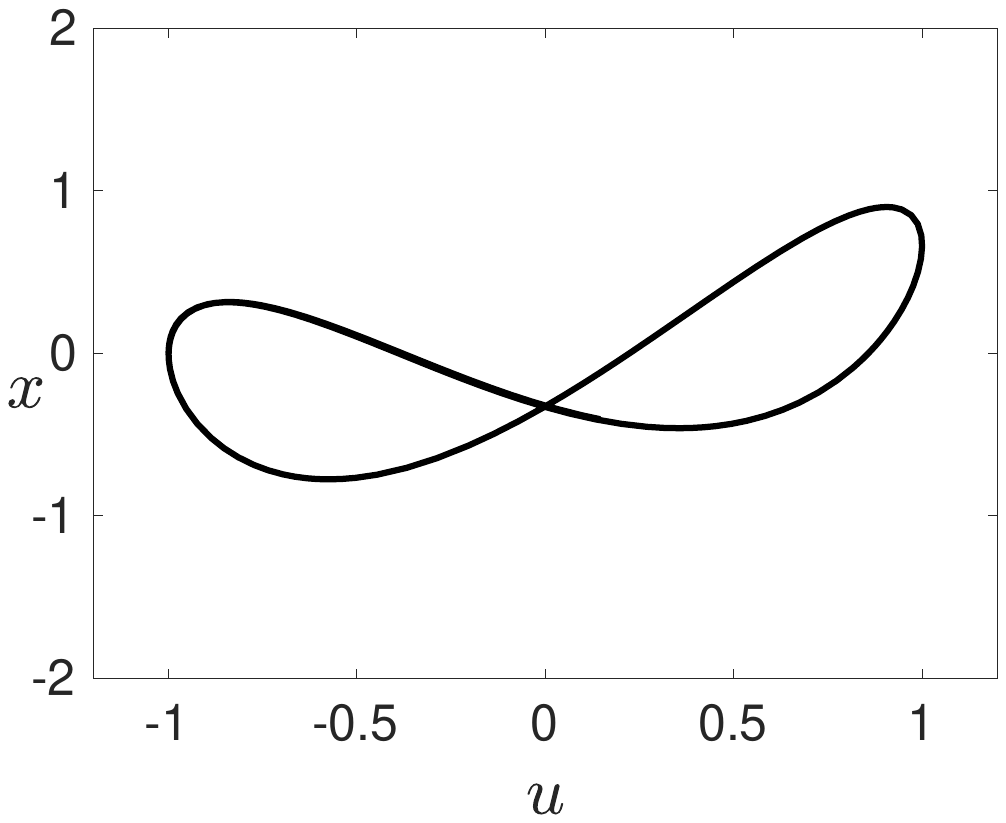}} \hspace{-0.5 cm}
\caption{\small  For a single differential equation with solution $x \in \mathbb R$ and periodic input $u\in \mathbb R$, displayed is the graph of $x$ as a function of $u$. The graph is that of a   closed curve that is not simple; that is, the loop formed crosses or ``pinches" itself. }
\label{figpinched}
\end{center}
\end{figure}

Consider for a single ODE with input $u$ and solution $x$, the curves formed by a periodic input $u$ versus $x.$ (If there is  a system of ODEs, some linear combination of the solutions $x_i$ is plotted against $u$.) In \cite{Bernstein2005} it is noted that the presence of a  nontrivial closed curve that {\it persists} as the frequency of the input approaches zero is characteristic of hysteresis. This nontrivial closed curve is called a{\it{ hysteresis loop}}.  
 If the loop degenerates such that the graph can be described by a function, this is a trivial closed curve. Figure~\ref{fig:comparehysteresis} contrasts a nontrivial closed curve and a trivial closed curve. Figure~\ref{figpinched} depicts an example of a pinched hysteresis loop  \cite{Wang2016}. These are closed curves that are not simple; that is, the loop crosses itself.  In \cite{Bernstein2011},  the authors study systems that exhibit hysteresis loops shaped like butterfly wings, some of which could also be classified as pinched loops. 
 It is important to note the presence of looping is not sufficient to conclude hysteresis; the loops must persist as the frequency of the input goes to zero.

 To establish hysteresis in the dynamical systems presented in this paper, the following procedure is applied.  First, the input $u$ is set to a constant and the equilibria of each system and their stability is determined.  Determining the existence of persistent looping is also useful. The hysteresis loops are generated using MATLAB and the stiff ODE solver ode23t  is used. In this paper, the periodic input $u(t)$ is typically of the form $\sin(\omega t);$ so, $\omega$ is the frequency of the input.

\begin{example} \exsprucebudworm Consider 
\begin{align}\label{exbudworm}
\dot x(t)& =|u(t)|x(t) \left(1 - \frac{x(t)}{q}\right) - \frac{(x(t))^2}{1 + (x(t))^2},
 \end{align}
where $x(t)$ is the population of the spruce budworm, $u(t)$ depends on the birth rate of the spruce budworms and the constant $q$ depends on the carrying capacity of the environment \cite{Ludwig1978}, \cite[Page 7-8]{Murray2002}. The carrying capacity is how much of the population the environment can support.

\begin{figure} [h] 
\centering
\includegraphics[width = .25\linewidth]{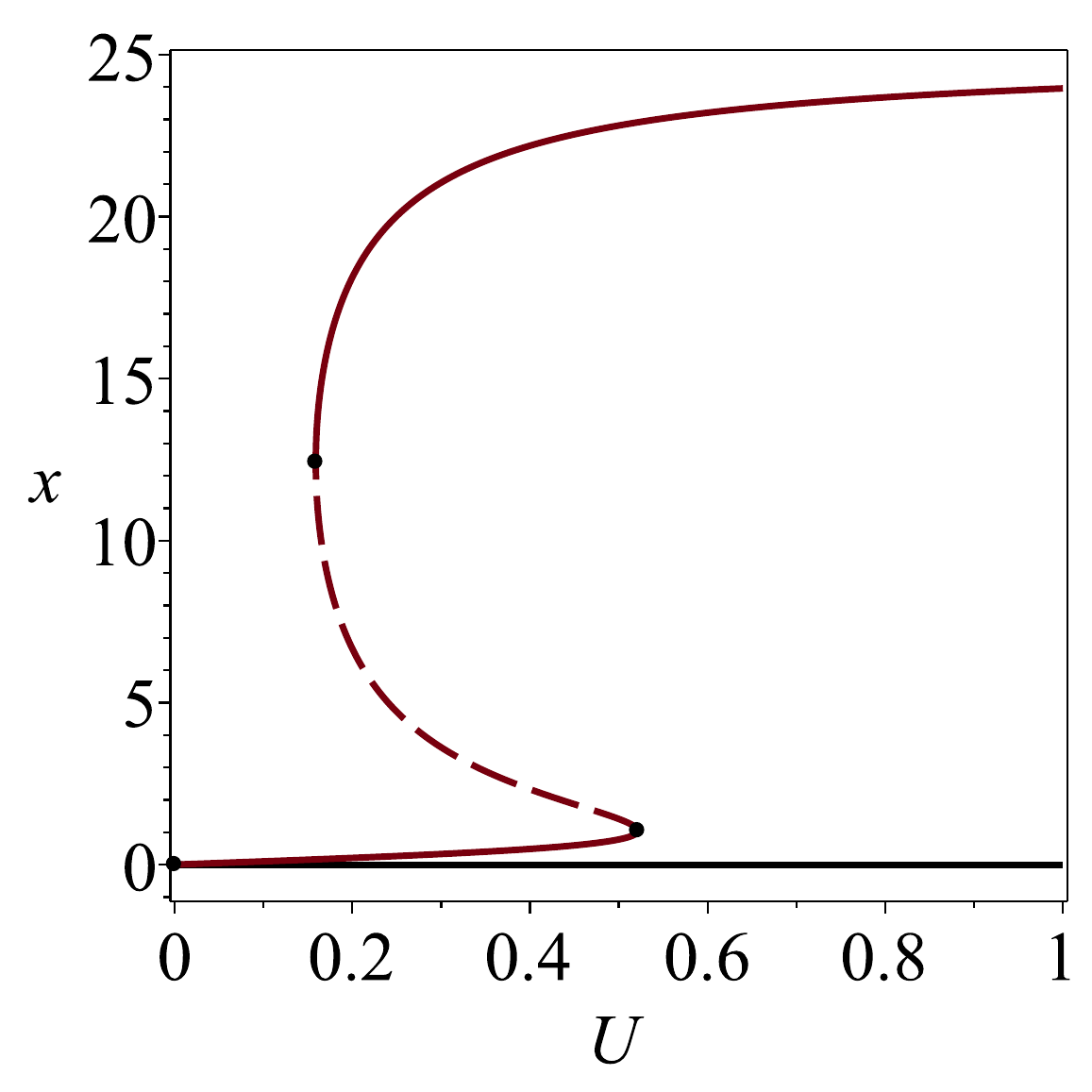}
\caption{\small \exsprucebudworm This diagram shows the equilibria of \eqref{exbudworm} and their stability, where the dashed line indicates an unstable portion of the graph. The change in stability occurs at $U = 0, 0.159, 0.521$. The black line indicates how regardless of the value of $U$, there is always an equilibrium point at 0.  \label{BudwormBifurcation} }
\end{figure}

As in Example~\ref{exmotivation},  various values of constant $u(t)$ affect the equilibria of \eqref{exbudworm} as shown in  Figure~\ref{BudwormBifurcation}. This figure indicates changes in the quantity and stability of equilibria occur when $u$ equals $0, 0.159$  and $0.521$. For instance, when $u(t)$ is set to be the constant input 0.52 and $q=25$, equation~\eqref{exbudworm} exhibits multiple stable equilibria; namely,
\begin{align*}
\bar a_1 &= 0,\\
\bar a_2 &= 0.972,\\
 \bar a_3 &= 1.123, \mbox{ and }\\
 \bar a_4 &=22.905
 \end{align*}
 as computed in Maple. As well, $\bar a_1, \bar a_3$ are unstable equilibria and $\bar a_2, \bar a_4$ are stable.  The phase portrait, constructed using Maple, in Figure~\ref{fig:BudwormPhase} supports this. 
\begin{figure}[h!]
\begin{center}
 \subfloat{\includegraphics[width = 0.35 \linewidth]{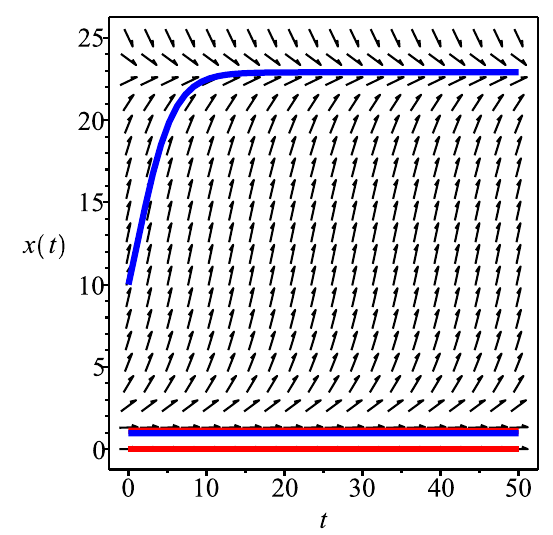}}\qquad
  \subfloat{\includegraphics[width = 0.35 \linewidth]{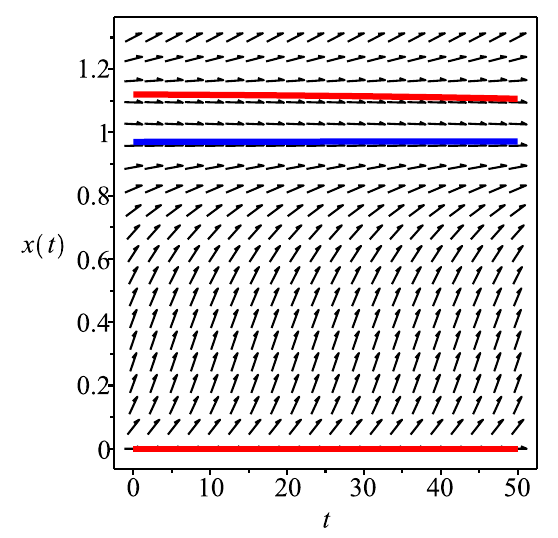}}
\caption{\small \exsprucebudworm  A phase portrait of \eqref{exbudworm} illustrating  $\bar a_1=0, \bar a_3=1.123$ are unstable equilibria and $\bar a_2= 0.972, \bar a_4=22.905$ are stable. The figure at the right is zoomed in because the equilibrium points are very close together near 0. The red curves indicate solutions near the unstable equilibria, while the blue curves indicate solutions near the stable equilibria.  \label{fig:BudwormPhase}}
\end{center}
\end{figure}

 The corresponding hysteresis loops of \eqref{exbudworm} are shown in Figure~\ref{fig:Budworm} and indicates persistent looping as $\omega$ goes to zero.  Figure~\ref{fig:BudwormZoomed} shows  the vertical ``jumps" of the hysteresis loops of \eqref{exbudworm} approaching closer to $0.159$  and $0.521$.

\begin{figure}[h!]
\begin{center}
\subfloat[]{\label{Ba}\includegraphics[width = 0.24 \linewidth, trim=0 175 0 175]{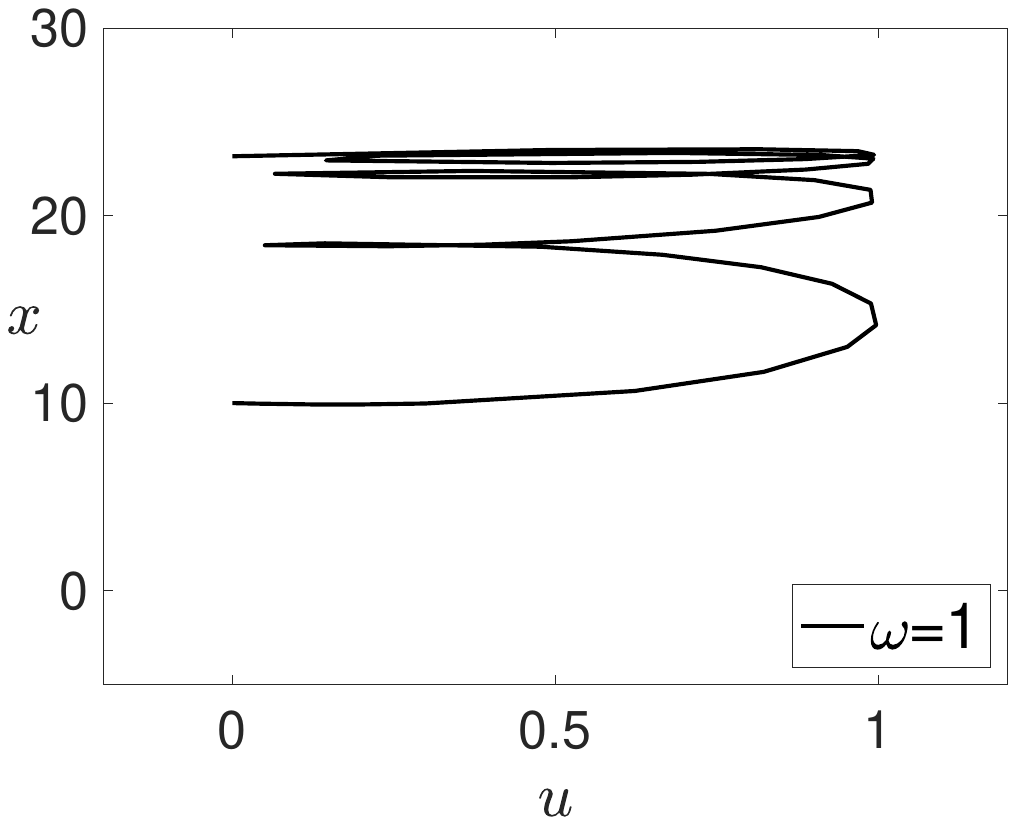}}
 \subfloat[]{\label{Bb}\includegraphics[width = 0.24 \linewidth, trim=0 175 0 175]{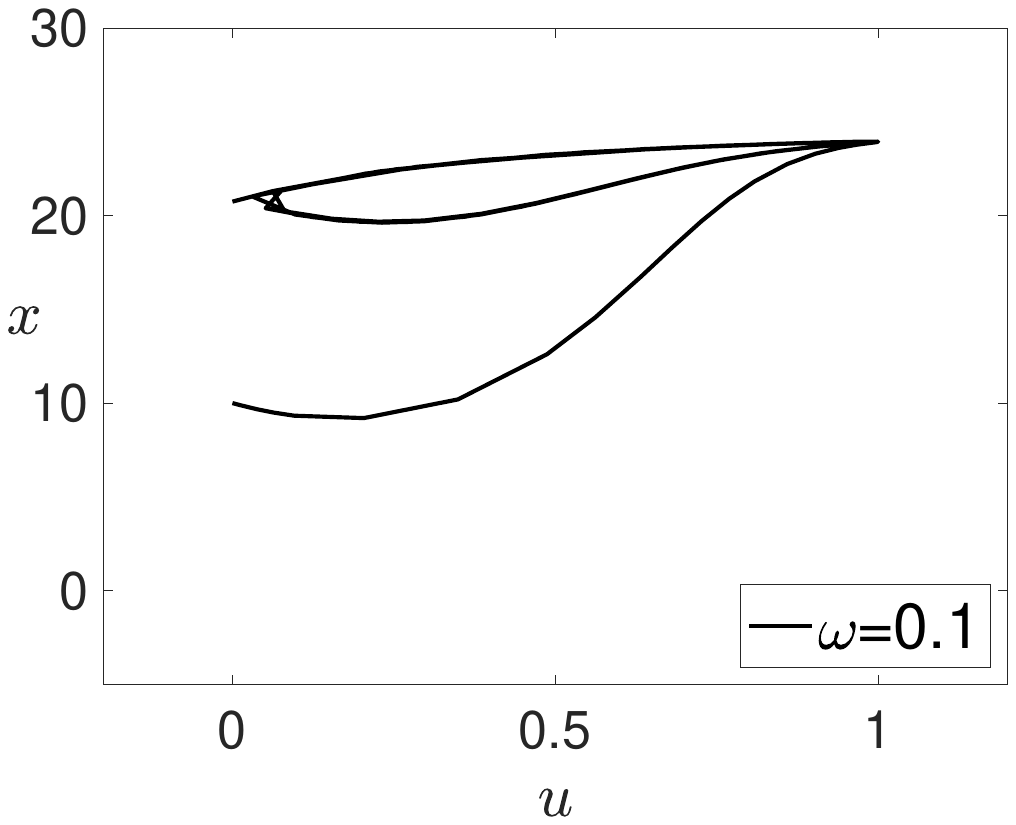}} 
 \subfloat[]{\label{Bc}\includegraphics[width = 0.24 \linewidth, trim=0 175 0 175]{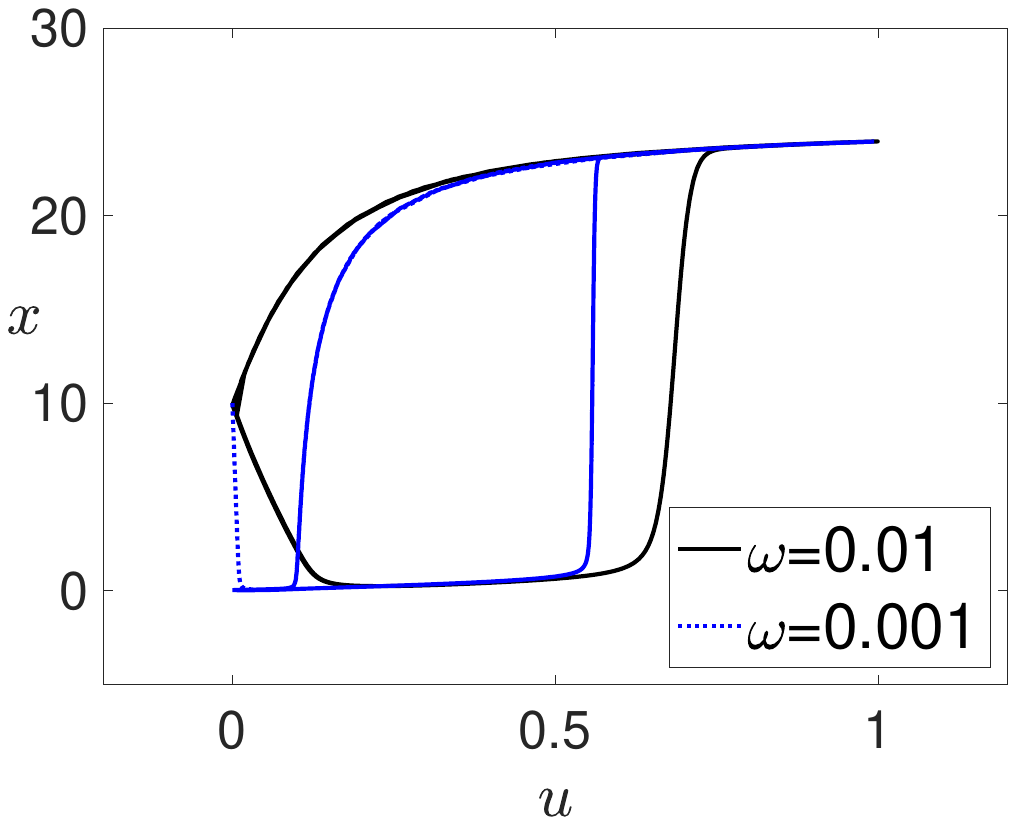}}
 \subfloat[]{\label{Bd}\includegraphics[width = 0.24 \linewidth, trim=0 175 0 175]{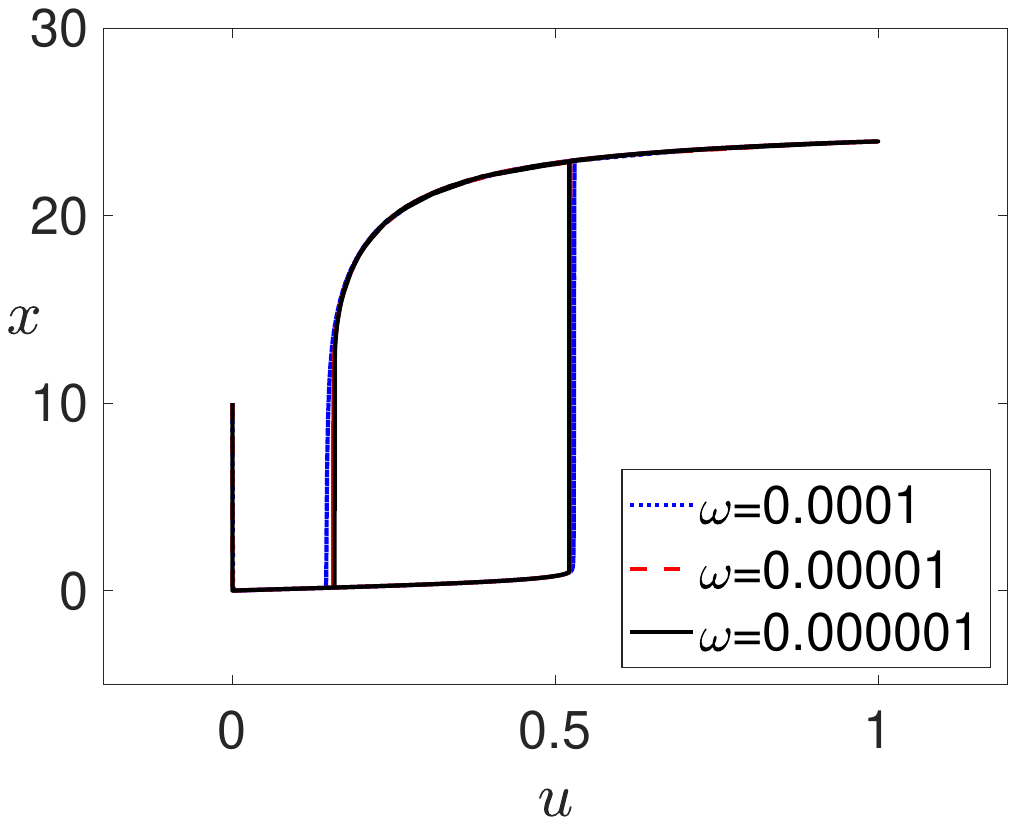}}
\caption{\small \exsprucebudworm  The graph of the solution, $x$, of \eqref{exbudworm} as a function of the input $u$. In this case, $u(t)=|\sin(\omega t)|$ for various values of $\omega$,  the  initial condition is $x(0)=10$ and $q=25$. Equation~\eqref{exbudworm} has two stable equilibria and two unstable equilibria.} 
\label{fig:Budworm}
\end{center}
\end{figure}

\begin{figure}[h!]
\begin{center}
 \subfloat{\includegraphics[width = .35 \linewidth,trim=0 175 0 175]{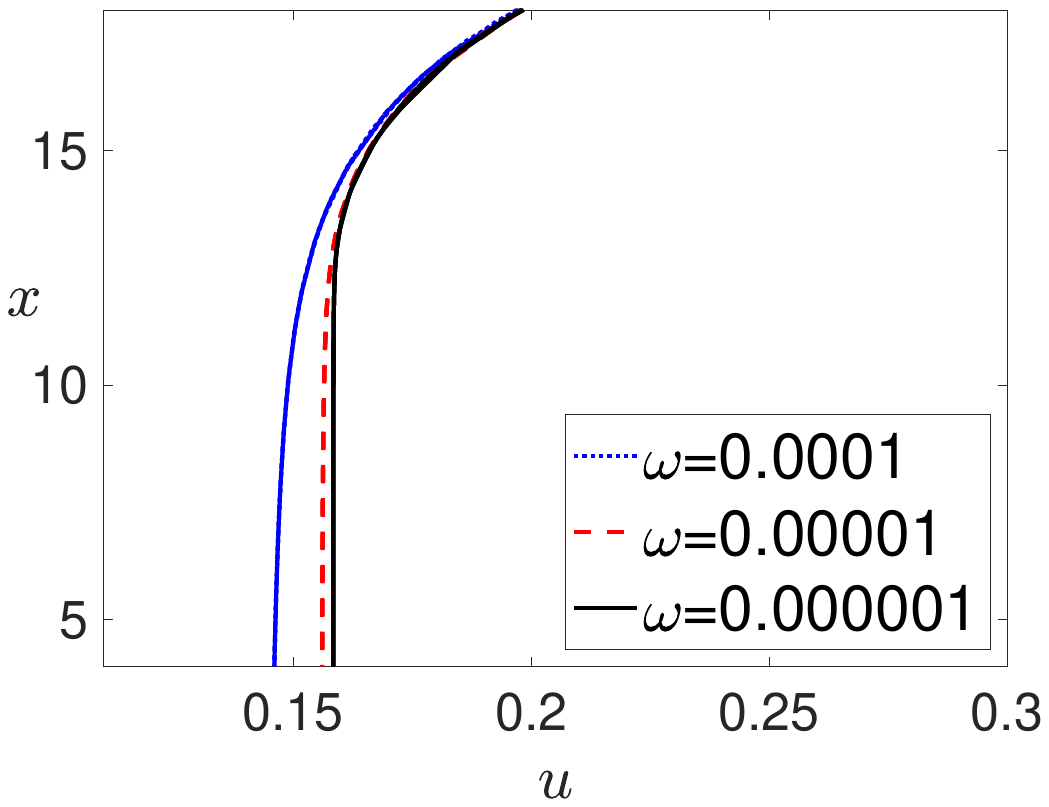}} 
 \subfloat{\includegraphics[width = .35 \linewidth,trim=0 175 0 175]{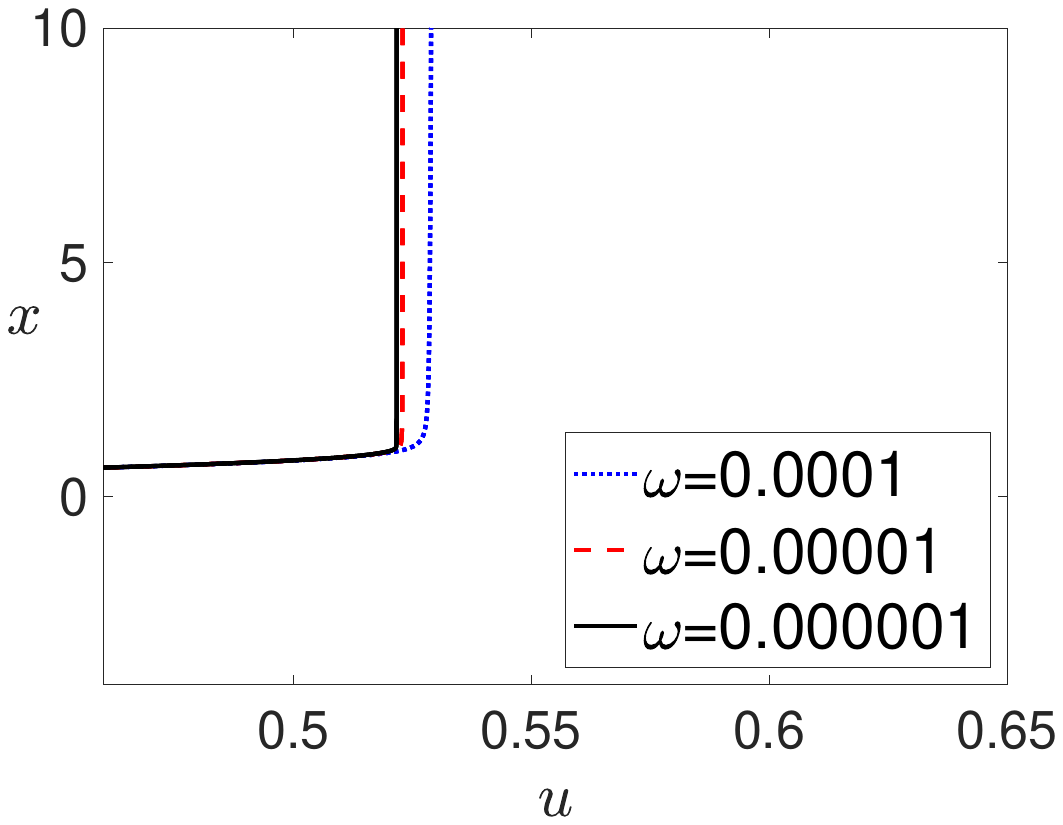}} 
\caption{\small \exsprucebudworm  A magnified image of Figure~\ref{Bd} shows the vertical jumps of the hysteresis loops converge to $0.159$  and $0.521$ as the frequency of the input $\omega$ approaches $0$.  \label{fig:BudwormZoomed}}
\end{center}
\end{figure}

\end{example}

\begin{example} \exflexiblebeam 
Consider 
\begin{align}\label{magbeam}
\ddot y(t)&= -\dot y(t)+0.5y(t)(1-(y(t))^2)+u(t).
\end{align}
This equation describes the behaviour of a flexible beam in a magnetic field \cite{Moon1979}. If the beam is composed of a ferromagnetic material, the magnetic field can bend the beam and $y(t)$ denotes the distance of the beam from its original position due to the magnetic field.  For an illustration of this, see Figure~\ref{figmagbeamsetup}.

\begin{figure}[h!]
\begin{center}
\begin{tikzpicture}
\draw[orange, fill](2,0)--(3,0)--(3,0.25)--(2,0.25);
\draw[orange, fill](5,0)--(6,0)--(6,0.25)--(5,0.25);
\draw[orange, ultra thick](3,0)--(4,0)node[anchor=north]{magnets};
\draw[black, ultra thick] (4,4)--(4,3)node[anchor=west]{flexible beam};
\draw[gray, ultra thick] (7,0)--(1,0)--(1,4)--(4,4);
\draw[gray, dashed, thick] (4,4)--(4,1);
\draw[gray,->](4,1.1)--(4.3,1.1)node[anchor=north]{$y$};
\draw[black, ultra thick] (4,4)--(4,3)--(4.4,1);
\end{tikzpicture}
\caption{\small \exflexiblebeam Magnets applied to a flexible ferromagnetic beam.}
\label{figmagbeamsetup}
\end{center}
\end{figure}

Defining $x_1=y(t)$ and $x_2=\dot y(t)$, and re-writing \eqref{magbeam} leads to
\begin{align*}
\dot x_1(t)&= x_2(t)\\
\dot x_2(t)&=-x_2(t)+0.5x_1(t)(1-(x_1(t))^2)+u(t).
\end{align*}
When $u(t)=0$, the equilibria are $(\pm 1,0)$, $(0,0)$,  and the eigenvalues of the corresponding Jacobian matrices indicate $(\pm 1,0)$ are both stable equilibria, while $(0,0)$ is an unstable equilibria \cite{Moon1979,Morris2011}. This echos the behaviour  in example~\ref{exmotivation}.  For $y$ the solution to  \eqref{magbeam}, Figure~\ref{figmagbeam} shows the behaviour of $y$ as a function of the input  $u(t)=\sin(\omega t)$  for different frequencies $\omega$, and persistent looping is observed as $\omega$ goes to zero. In addition, the shape of the hysteresis loops are similar to those of example~\ref{exmotivation}, see Figure~\ref{fig:Perko5}. This is likely due to both examples having the same type of equilibria; namely, one unstable equilibria sandwiched between two stable equilibria, and the two stable equilibria are equidistant from the unstable equilibria.

\begin{figure}[h!]
\begin{center}
  \subfloat[]{\includegraphics[width = .24\linewidth, trim=0 175 0 175]{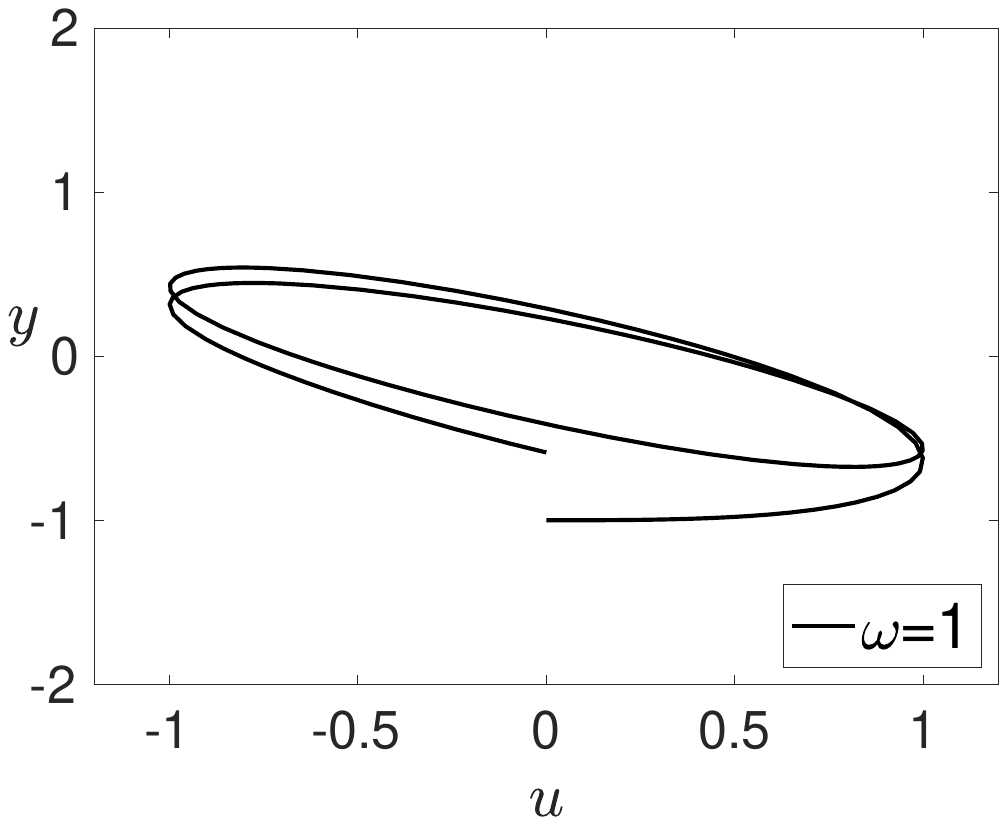}}
 \subfloat[]{\includegraphics[width = .24 \linewidth, trim=0 175 0 175]{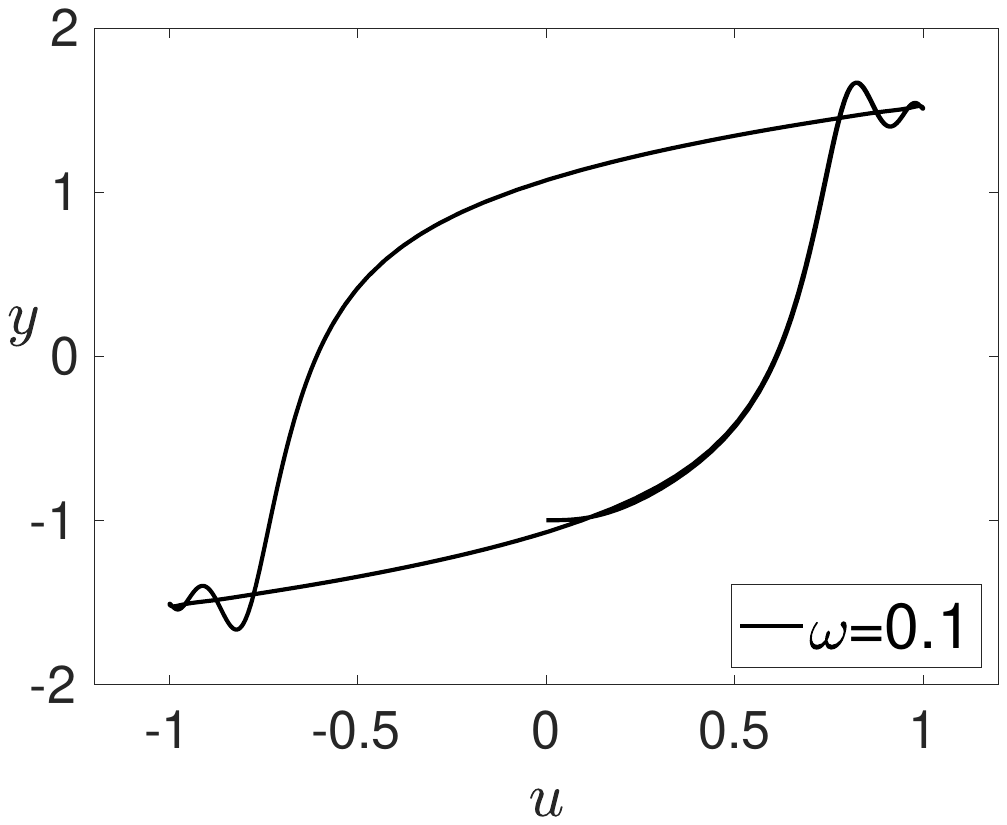}} 
   \subfloat[]{\includegraphics[width = .24 \linewidth, trim=0 175 0 175]{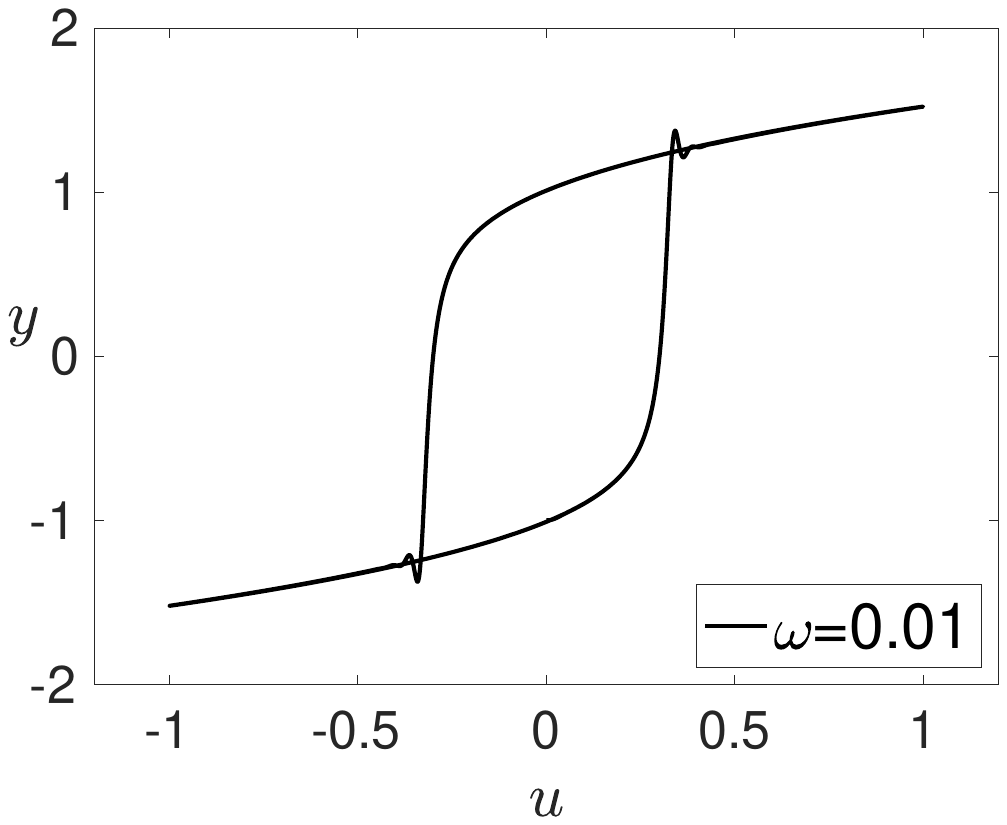}} 
\subfloat[]{\includegraphics[width = .24 \linewidth, trim=0 175 0 175]{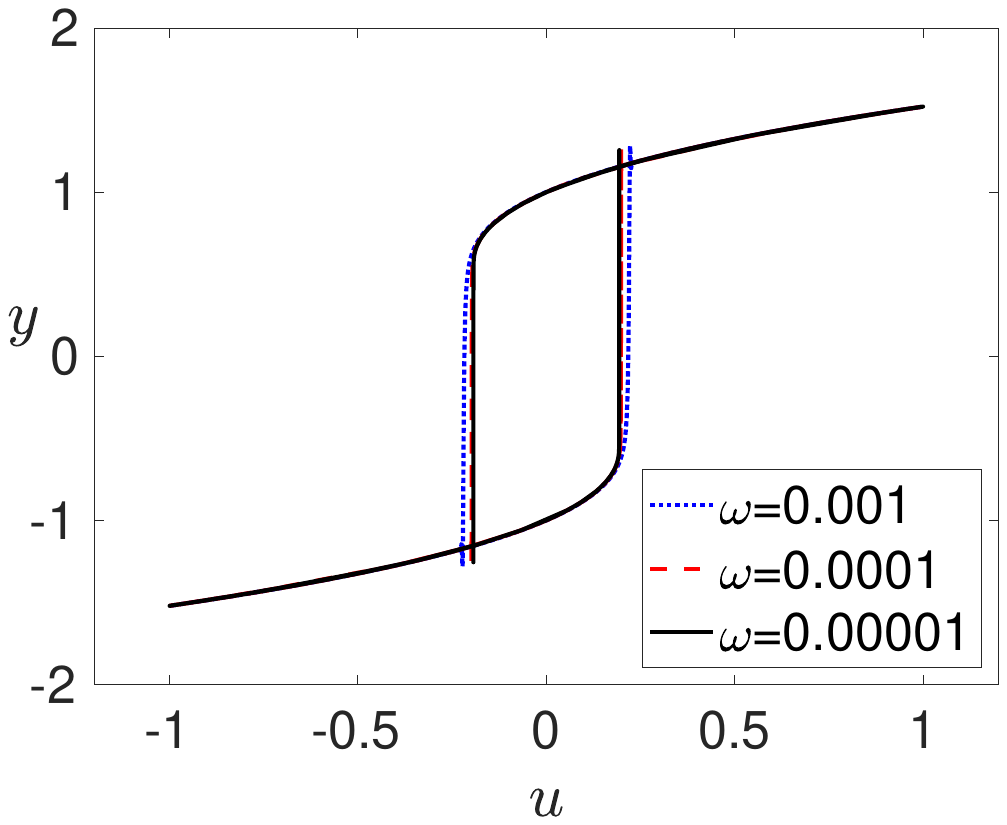}}
\caption{\small \exflexiblebeam The solution $y$ of \eqref{magbeam} as a function of $u(t)=\sin(\omega t)$ is shown for various values of $\omega$. Persistent looping occurs and as $\omega$ approaches zero. The initial conditions are $y(0)=-1$ and $y(0)=0$. Equation~\eqref{magbeam} has two stable equilibria and one unstable equilibrium point between the two stable equilibria.}
\label{figmagbeam}
\end{center}
\end{figure}

\end{example}

\begin{example} \expendulum 
\begin{figure}[h!]
\begin{center}
\begin{tikzpicture}
\draw(4.4,2.8) node[anchor=west]{$l$};
\draw[gray, dashed, thick] (4,4)--(4,0)node[anchor=north]{$y(0)=0$};
\draw[ultra thick] (4,4)--(5,0.7)node[anchor=west]{$m$};
\filldraw (5,0.7) circle (3pt);
\filldraw[gray] (4,4) circle (1pt);
\draw[gray,thick, ->] (4,1) arc (90:110:-2.45);
\end{tikzpicture}
\caption{\small  \expendulum Depicted is the dynamics of a simple pendulum. It has length $l$ and a mass $m$ at its end. The angle of the pendulum from its vertical rest position is $y$ and at its initial condition, it is at rest; that is $y(0)=0$.  }
\label{figpendulum}
\end{center}
\end{figure}
The equation of motion of a simple pendulum with friction and forcing term $u(t)$ is 
\begin{subequations}\label{eqsecondorder}
\begin{align}
\ddot y(t) &=-\frac{k}{m} \dot y(t) -\frac{g}{l}\sin(y(t)) + u(t)
\end{align}
\end{subequations}
where $y(t)$ is the angle of the pendulum from its vertical rest position, $m$ is the mass hanging from the pendulum of length $l$, and $g$ is acceleration due to gravity. The frictional force is proportional to the speed of the mass and $k$ is the coefficient of friction.  When $u(t)=0$, this equation has two equilibria, $(0,0)$ and $(\pi,0)$, and only $(0,0)$ is stable. Hence, this system cannot exhibit hysteresis.  Figure~\ref{fig:Pendulum} demonstrates persistent looping is absent. 

\begin{figure}[h!]
\begin{center}
\includegraphics[width = 0.24\linewidth , trim=0 175 0 175]{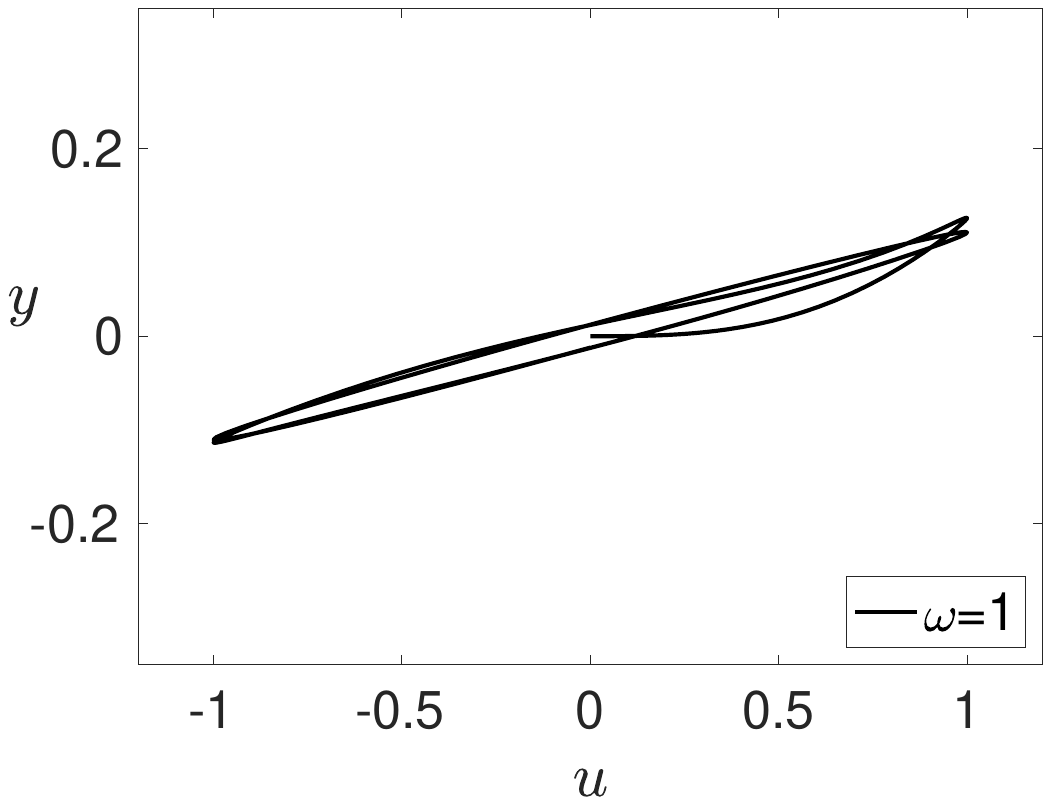}
\includegraphics[width = 0.24\linewidth , trim=0 175 0 175]{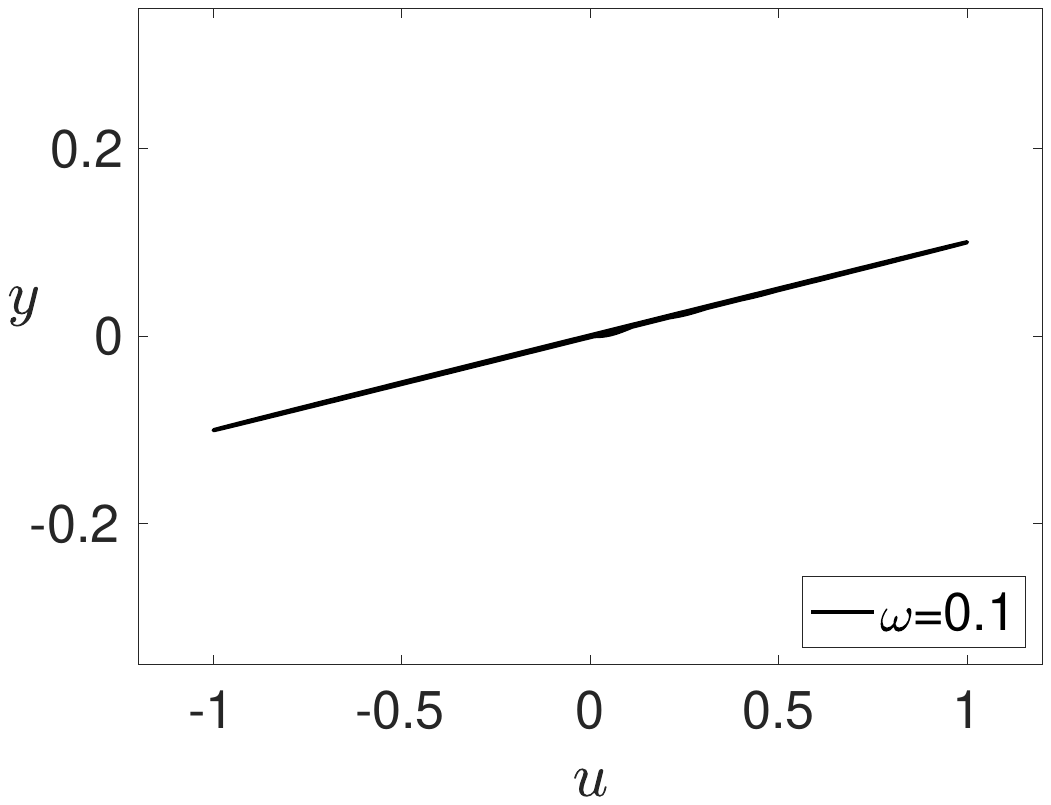}
\includegraphics[width = 0.24\linewidth , trim=0 175 0 175]{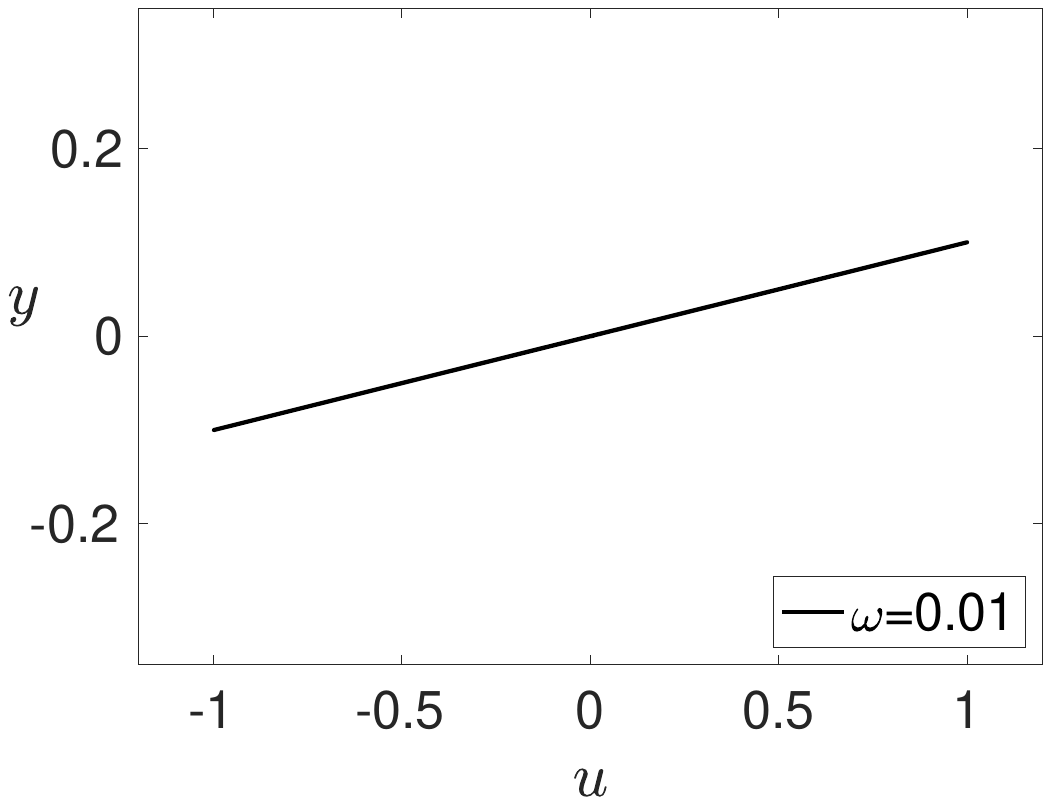}
\includegraphics[width = 0.24\linewidth , trim=0 175 0 175]{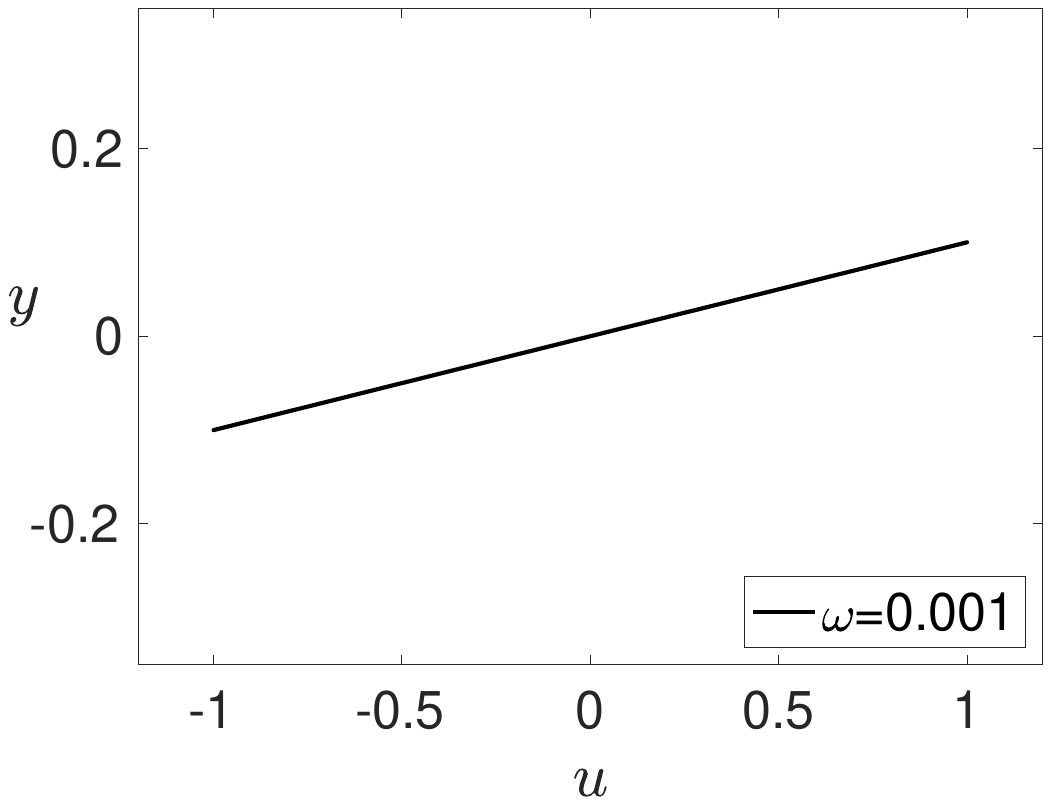}
\caption{\small \expendulum Graph of the solution $y$ of \eqref{eqsecondorder} as a function of the input $u$. In this case,  $u(t)=\sin(\omega t)$ for various values of $\omega$, and the initial condition is $y(0)=0, \dot y(0) = 0$. The parameters in \eqref{eqsecondorder}  are $k=m=1$ and $g=10$. This system only has two equilibria and of the two, only one is stable. Persistent looping does not occur. }\label{fig:Pendulum}
\end{center}
\end{figure}
\end{example}

The hysteresis examples thus far are typical in that they have isolated stable equilibria and are described by nonlinear dynamics.  In the following  example,   there is a continuum of stable equilibria and the nonlinearity in the differential equation is only in the input.

\begin{example} \label{exDuhemModel} \exferro \cite{Coleman1986,Bernstein2005}. The Duhem Model  describes ferromagnetism and its corresponding equation is
\begin{align}\label{eqDuhem}
\dot x(t) &= \alpha |\dot u(t)| (\beta u(t) - x(t) ) + \gamma \dot u(t).
\end{align}
Notice equation \eqref{eqDuhem} is linear with respect to $x(t)$. The input $u(t)$ is the magnetic field strength and the solution to \eqref{eqDuhem} is the magnetic flux density \cite{Coleman1986,Bernstein2005}. The parameters $ \beta, \gamma$ are positive constants. The equilibria of \eqref{eqDuhem} are determined by setting the input $u(t)$ to be constant, meaning $\dot u(t)=0$. This leads to $\dot x(t) = 0$ for all $t$.   If $\alpha <0$ then the continuum of equilibria is unstable and hence  \eqref{eqDuhem} cannot exhibit hysteresis.  Figure \ref{fig:DuhemN} displays no looping behaviour for \eqref{eqDuhem} for a specific instance of $\alpha<0$. On the other hand, Figure~\ref{fig:Duhem} depicts a scenario when $\alpha>0$ and looping behaviour does occur. Furthermore, notice the shape of the hysteresis loops do not change despite $\omega$ changing. 


\begin{figure}[h!]
\begin{center}
\includegraphics[width = .24\linewidth, trim=0 175 0 175]{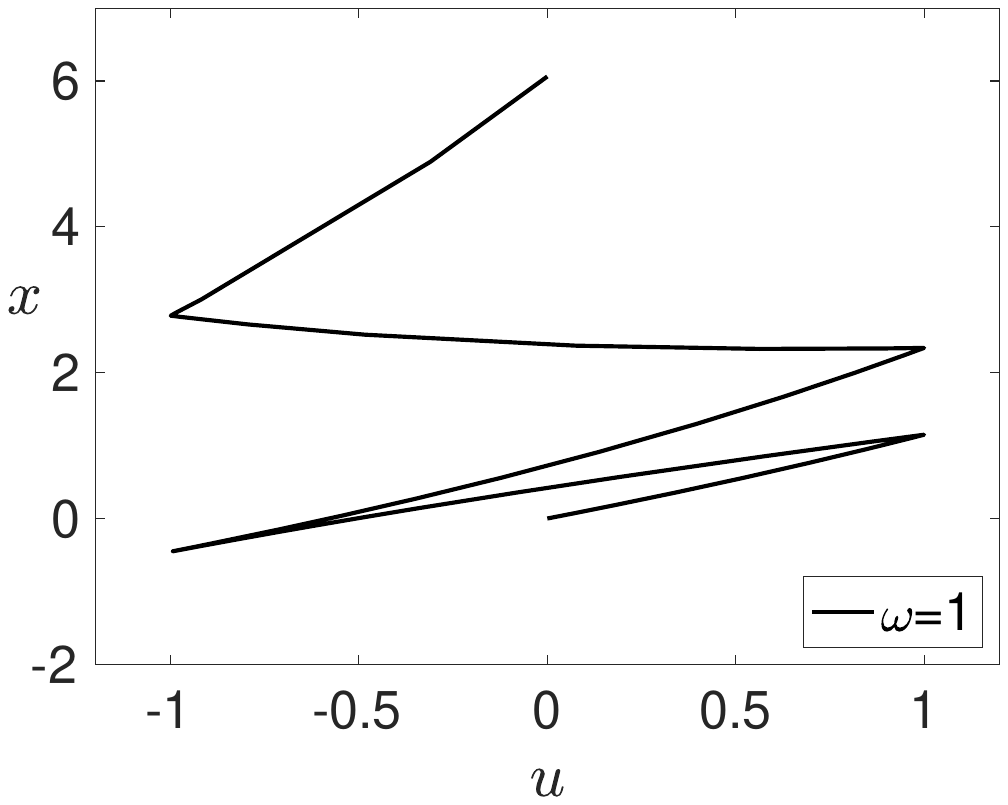}
\includegraphics[width = .24\linewidth, trim=0 175 0 175]{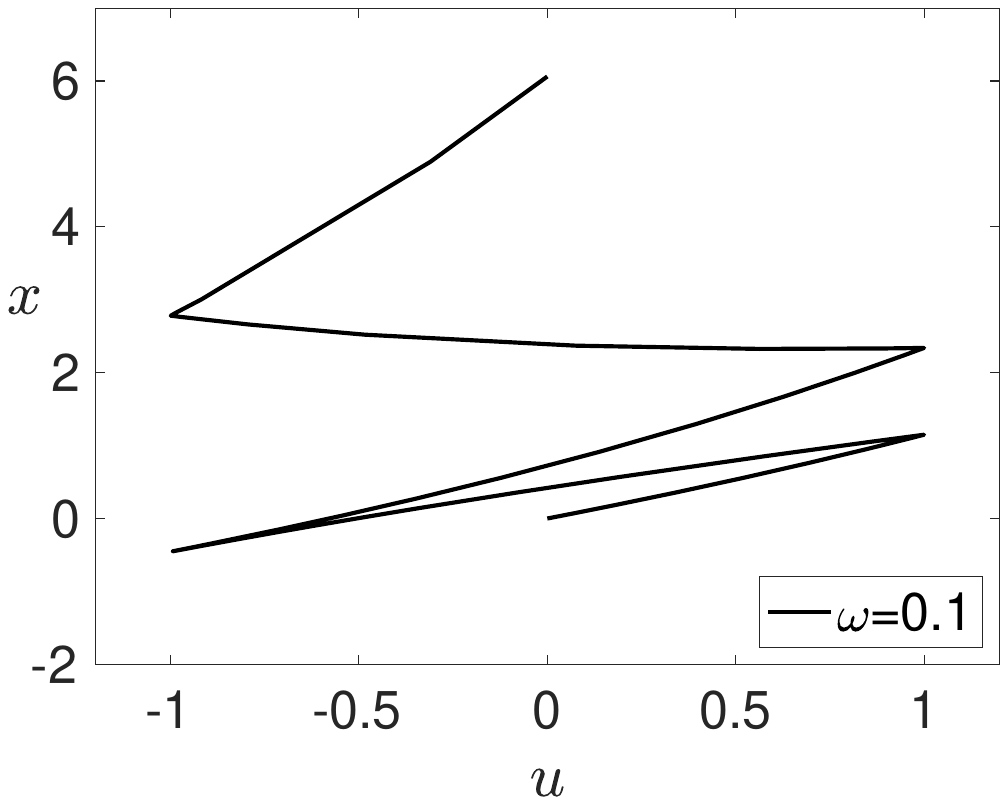}
\includegraphics[width = .24\linewidth, trim=0 175 0 175]{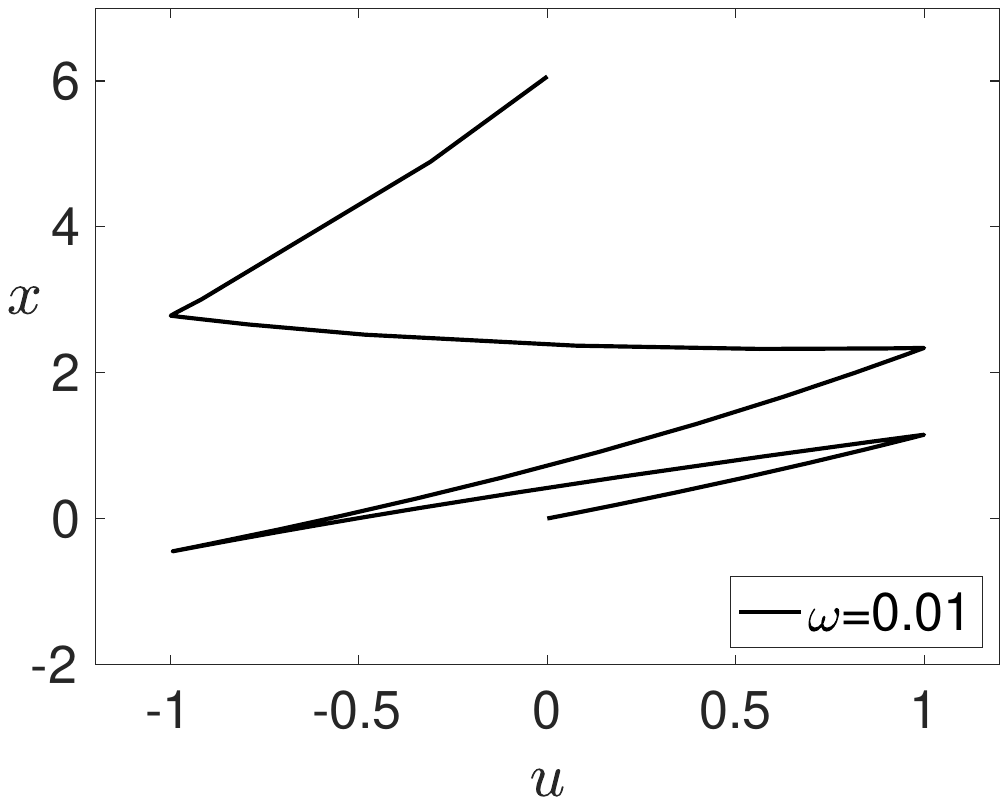}
\includegraphics[width = .24\linewidth, trim=0 175 0 175]{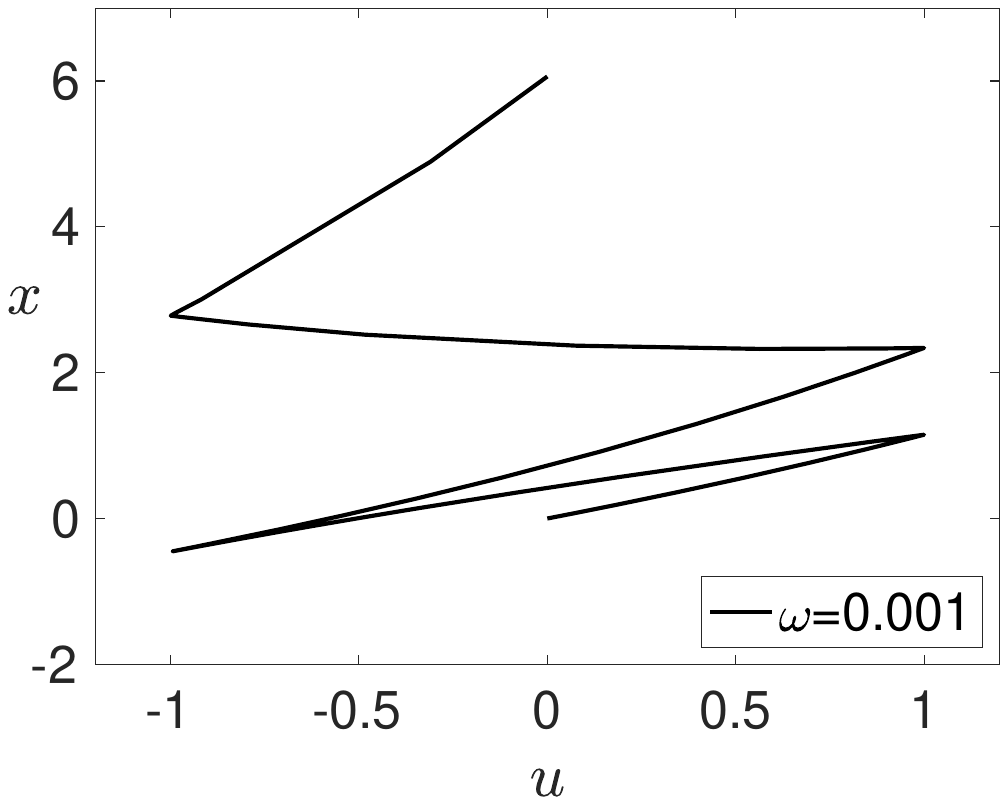} 
\caption{\small  \exferro  Graph of  the solution $x$ of  \eqref{eqDuhem}  as a function of the input $u$ with $\alpha=-1, \beta = 0.5, \gamma = 1, x(0) = 0$ and $u(t) = \sin(\omega t)$, where $\omega$ is the frequency of the input. No closed curves appear, and therefore the system does not exhibit hysteresis in the case of negative $\alpha$.}
\label{fig:DuhemN}
\end{center}
\end{figure}

\begin{figure}[h!]
\begin{center}
\includegraphics[width = .24\linewidth, trim=0 175 0 175]{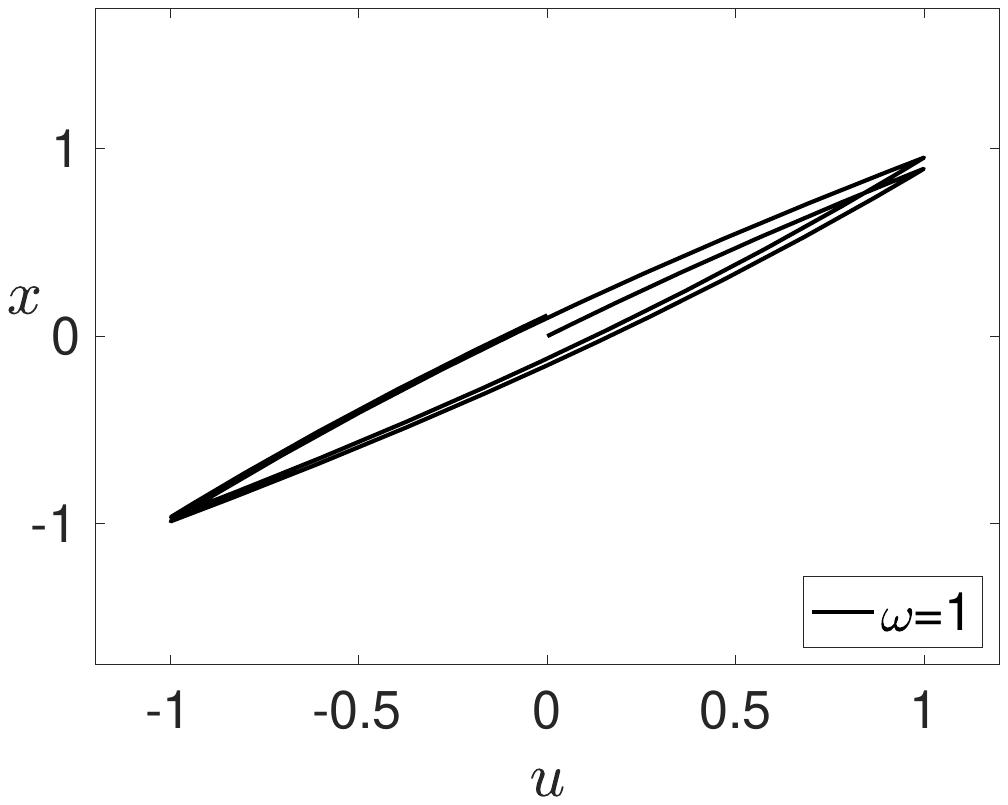}
\includegraphics[width = .24\linewidth, trim=0 175 0 175]{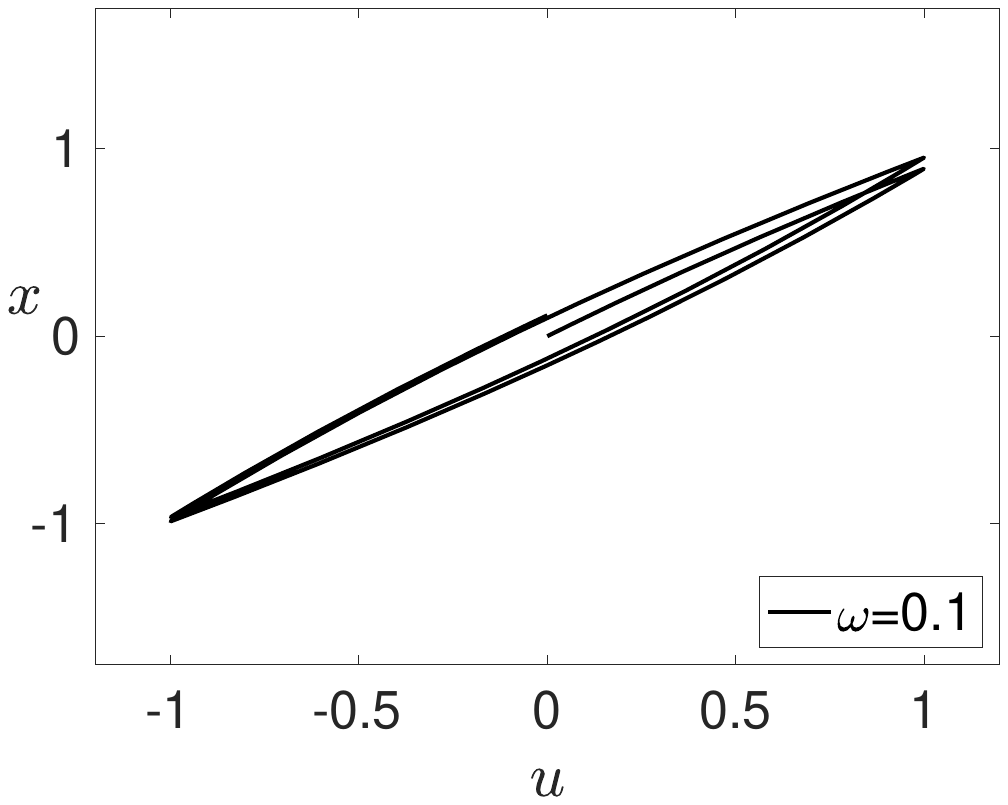}
\includegraphics[width = .24\linewidth, trim=0 175 0 175]{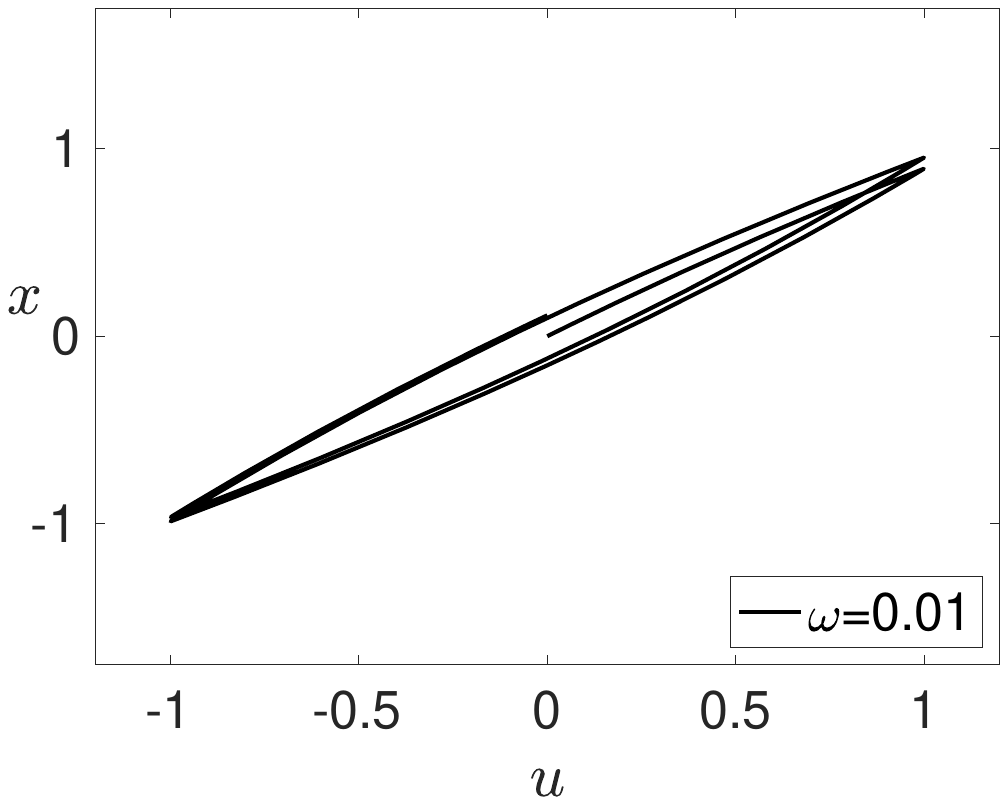}
\includegraphics[width = .24\linewidth, trim=0 175 0 175]{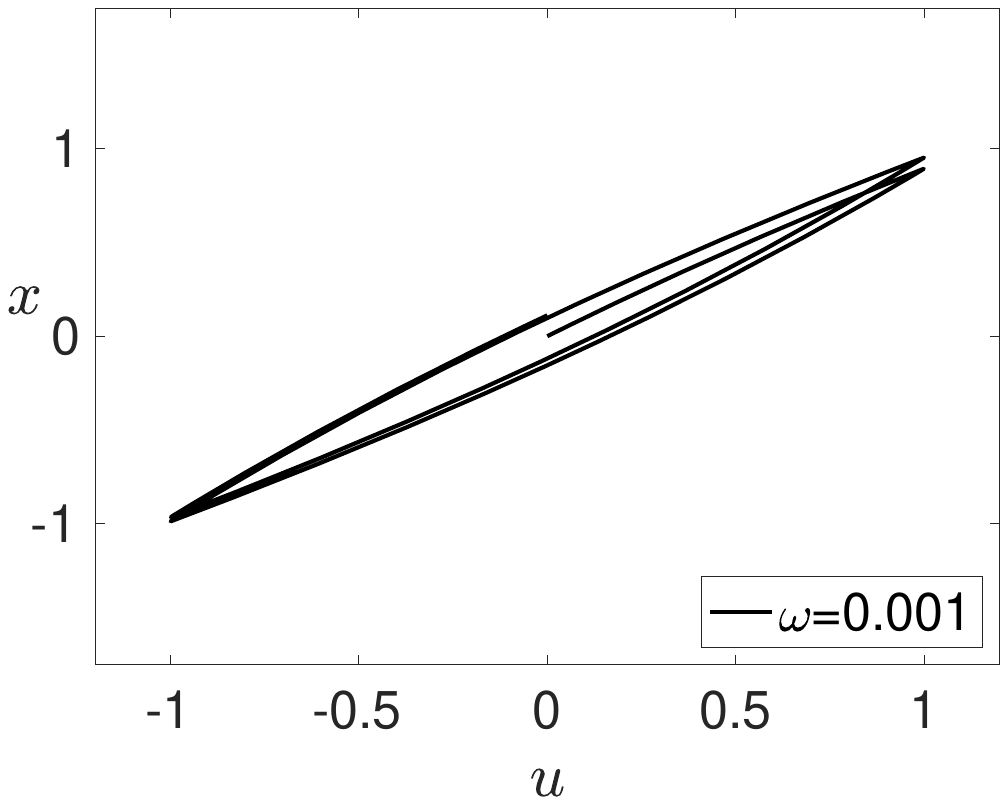}
\caption{\small  \exferro  Graph of the solution $x$ of  \eqref{eqDuhem}  as a function of the input $u$ with $\alpha=0.5, \beta = 0.5, \gamma = 1,$ $x(0) = 0$ and $u(t) = \sin(\omega t)$, where $\omega$ is the frequency of the input.  Persistent looping behaviour is demonstrated in the case of positive $\alpha$.}
\label{fig:Duhem}
\end{center}
\end{figure}
\end{example}

\section{Discussion}
Hysteresis has been analyzed here from a  dynamical systems viewpoint and illustrated with low-order differential equations. It is  a useful additional concept to a course covering nonlinear differential equations. 
The examples of ODEs exhibiting hysteresis  presented can be used as a basis for student exercises. Additional exercise problems are provided in the appendix. The exercise problems are a resource for assessing student understanding of hysteresis and related concepts such as stability. A number of the exercise problems use software like Maple and MATLAB, which re-enforce course learning objectives in  understanding dynamical systems and  computational expertise.  Lyapunov functions are important in stability analysis and can be used in exploring deeper connections between hysteresis and stability. Discussion and several examples can be found in  \cite{Morris2011}. 

Another point for investigation is to look for connections between the shape of the hysteresis loops and the system equilibria.  Operators are  useful tools for considering the shape of hysteresis loops \cite{Brokate1996,Visintin2005}. For example, the operator
\[
F(u)=\left\{\begin{array}{rr} 
-1, & -5<u<5\\
1, & 0<u<10
\end{array}\right\}
\]
can be used to describe the  hysteresis loop in  Figure~\ref{figrelay1}. This style of hysteresis operator  has been used to reduce frequent `on' and `off' switching from devices like thermostats \cite{Morris2011}. However,  when the input is between 0 and 5, it is not clear whether $x$ is 1 or -1, unless the previous value of $x$ is known.  For this reason, hysteresis is often said to have a memory, and that  the memory effect in hysteresis is related to its looping behaviour.  In a dynamical system, the state is the memory of the system and knowledge of $x(t_0)$ and $u(t)$ for $t\geq t_0$  is sufficient to determine  values of $x(t)$ for $t \geq t_0.$  The operator  approach does not explicitly indicate  the intrinsic properties of the underlying dynamical system. It is useful for complex systems like smart materials; see for example, \cite{Smith2005}. 
 \begin{figure}[h!]
\begin{center}
\includegraphics[width = 4.1 cm, trim=1cm 6.8cm 1cm 5cm , clip]{relay.pdf}
\caption{\small  Relay}
\label{figrelay1}
\end{center}
\end{figure}

\begin{figure}[h!]
\begin{center}
\vspace{-1.6cm} \subfloat{\includegraphics[width = 4.3 cm, trim=1cm 6.8cm 1cm 0cm , clip]{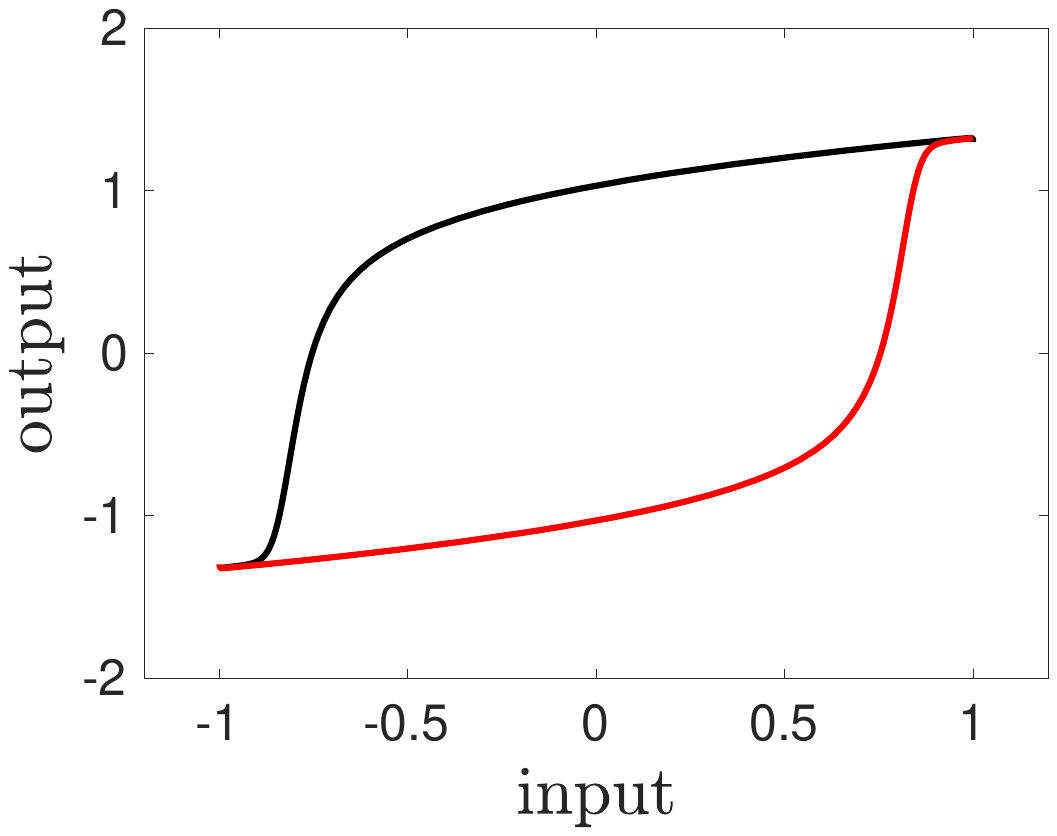}} \hspace{-0.5 cm}
\caption{\small  For the example of freezing and thawing, regard heat flow as input and the concentration of liquid/ice as output.  The freezing process  is not necessarily identical to the reverse process of thawing.  Consequently, if the black curve models the freezing process, and the red curve models the thawing process, then the dynamics of freezing and thawing forms a loop.}
\label{figpath}
\end{center}
\end{figure}

Hysteresis also occurs in partial differential equations (PDEs). The same considerations apply to hysteresis in PDEs as in the ODE case, such as multiple  stable equilibria and persistent looping behaviour. However, the analysis is more complex.   
The dynamics of freezing and thawing is described by a type of Boussinesq equation \cite{Alimov1998}. The freezing process is not necessarily identical to the reverse process of thawing and hence this creates a looping behaviour as depicted in Figure~\ref{figpath}.
The Landau-Lifshitz equation models the magnetic behaviour in ferromagnetic objects \cite{Guo2008,Landau1935}. It is a coupled system of  nonlinear PDEs  and is known to have multiple stable equilibrium points \cite[Theorem~6.1.1]{Guo2008} \cite{Chow2015}.   Persistent looping behaviour has been shown in \cite{Chow2014_ACC}.  The system of ODEs and PDEs describing  lithium ion batteries has multiple equilibria and simulations  also indicates persistent looping in the current-voltage relationship; see for example  \cite{Afshar2016}. Hysteresis also appears in electromagnetism, which is described by a special case of Maxwell equations \cite{Eleuteri2007}, and in population dynamics \cite{Aiki2005}.  

\appendix
\section{Exercise Problems}
In this appendix, exercise problems that complement the preceding discussion on hysteresis in ODEs is provided.

\begin{exercise} 
Consider the following set of isolated equilibria of \eqref{FONODE}
\[
S = \{\bar x_1, \bar x_2, ..., \bar x_n\}
\] 
for $n$ being some natural number. Show equation~\eqref{FONODE} for scalar $f$ cannot exhibit hysteresis when $\frac{df}{dx}(\bar x_j)>0$, where $\bar x_j \in S$ for all $j \in [0,n]$. 
\end{exercise}
\begin{exercise} 
Consider an equation of the form 
\begin{subequations}~\label{SONODE1}
\begin{align}
\ddot y(t)&= f(y(t),\dot y(t)),\\
y(0) &= y_0, \, \dot y(0)= y_1,
\end{align}
\end{subequations} 
where $f$ is a continuous function of $t\geq0$ and $y_0, y_1$ are constants. Assume $f$ and $y$ are scalar.  Rewrite \eqref{SONODE1} into a system of first order equations by setting $x_1(t)  = y(t)$ and $x_2 = \dot y(t)$. 
Let $n$ be some natural number, and consider a set of isolated equilibria  $S = \{\bar X_1, \bar X_2, ... \bar X_n \}$, where $\bar X_j =(x_{1j},x_{2j} )\in S$ for $j\in [0,n]$ is an equilibrium point of the first order system of \eqref{SONODE1}.  Show equation \eqref{SONODE1} cannot exhibit hysteresis if for all $j\in [0,n]$, the equilibrium $\bar X_j \in  S$ satisfies $\frac{\di f}{\di x_1}(\bar X_j)> 0 $ or $\frac{\di f}{\di x_2} (\bar X_j) >0.$ 
\end{exercise} 

\begin{exercise} \label{exercisecode}
Rewrite $\ddot y(t)+c\dot y(t)+k(y(t)-(y(t))^3)=u(t), y(0)=1, \dot y(0)=0$ into a system of first-order equations, and then use the following MATLAB function to construct the corresponding hysteresis loops for $u(t)=\sin(\omega t)$ for various values of $\omega$. 

\begin{verbatim}
function [t,x,input]=SpringDynamics(tspan,x0,a)

%tspan is the interval of time
%x0 is the initial condition, which has two components
%a is the frequency of the input

k=-1;
c=15;b=k;

[t,x]=ode23t(@nonlinearspring,tspan,x0);
input=u(t);
  
function dx=nonlinearspring(t,x)
dx1=x(2);
dx2=-c*x(2)-k*x(1)+b*x(1)^3+u(t);
dx=[dx1;dx2];
end;

function f=u(t)
f=sin(a*t);
end;
end
\end{verbatim}\\
\end{exercise} 

\begin{exercise} 
Consider 
\begin{subequations}~\label{FOLODE}
\begin{align}
\dot x(t) &= a(t)x(t) + b(t),\\
x(0) &= x_0,
\end{align}
\end{subequations}
where $a(t)$ is a continuous function of $t$  with $t\geq0$, $x_0$ is an arbitrary constant, and $ b(t)$ is some continuous function of the input.  Assume $x$, $a$ and $b$ are scalars.

\noindent (a) When $a(t)$ is a nonzero constant and $b(t)$ is a constant, show \eqref{FOLODE} cannot exhibit hysteresis.\\
(b) Suppose $a(t)=0$ and $b(t)=\sin(\omega t)$. Demonstrate the absence of hysteresis by modifying the MATLAB function in exercise~\ref{exercisecode}.\\
(c) For $a(t)=0$ and constant input, show \eqref{FOLODE} has a continuum of equilibria.\\
(d) Let $u(t)=\sin(\omega t)$, $a(t) = 0, b(t)= |\dot u(t)|u(t)$ and $x_0 = 0.5$.  Modify the MATLAB function in exercise~\ref{exercisecode} to demonstrate the presence of hysteresis for various values of $\omega$.

 Suppose the first part of definition~\ref{defMorris} is relaxed. In particular, instead of an equilibrium which is necessarily constant for all time;  consider {\it steady state solutions,} which only require convergence to a constant solution after some $t$. An example of a steady state solution could be $e^{-t}\cos(t)$ since this goes to zero as $t$ approaches infinity.   The following explore whether steady state solutions can lead to the existence of hysteresis: modify the MATLAB function in exercise~\ref{exercisecode} to demonstrate the presence of hysteresis for various values of $\omega$ in  \eqref{FOLODE} if  $b(t)= |\dot u(t)|u(t)$ with  $u(t)=\sin(\omega t)$, initial condition $x_0=0.5$ and \\
\noindent (e)  $a(t) = -e^{-t}$,\\
(f)  $a(t) =-\frac{1}{t+1}$,\\
(g) $a(t) =-\frac{1}{t+1}+1$, and\\
(h) $a(t) =-t$.
\end{exercise} 

\begin{exercise} For $x, u\in \mathbb R$,
\begin{equation}\label{egnonlinearcontinuum}
\dot x(t) = |\dot u(t)|(x(t) - (x(t))^3 + u(t)),
\end{equation}
show for any constant $u(t)$ that any constant is an equilibrium point and the equilibrium points  are all stable. Then for $u(t)=\sin(\omega t)$, modify the MATLAB function in exercise~\ref{exercisecode} to demonstrate the presence of hysteresis for various values of $\omega$. The resulting hysteresis loops will have the same shape despite the frequency of the periodic input changing, and this is known as  rate-independent hysteresis; otherwise, the hysteresis is said to be rate-dependent \cite{Bernstein2005}.  Some authors define rate independence as necessary for hysteresis to occur \cite{Brokate1996,Grinfield2009,Mielke2015,Visintin2005}. 
\end{exercise} 

\begin{exercise}   For $y, u\in \mathbb R$, consider
\begin{subequations}\label{eqDuhem2D} 
\begin{align}
\ddot y(t) &+p(t)\dot y(t) +\alpha|\dot u(t) | y(t)=\alpha\beta |\dot u(t)|u(t)+~\gamma \dot u(t),\\
y(0) &= y_0, \dot y(0) = y_1.
\end{align}
\end{subequations}
(a) Rewrite \eqref{eqDuhem2D} into a system of first-order equations.\\
(b) Show for any constant $u(t)$ that any constant is an equilibrium point.\\
(c) For $u(t)=\sin(\omega t)$, modify the MATLAB function in exercise~\ref{exercisecode} to demonstrate the presence of hysteresis for various values of $\omega$ with\\
(i) $p(t)=0.5$, $\beta=\alpha=1$ and $\gamma=0$,\\
(ii) $p(t)=0.2$, and $\beta=\alpha=\gamma=1$, \\
(iii) $p(t)=e^{-t}+1$, $\beta=\alpha=1$ and $\gamma=0$. \\
Comment and compare the shape of the hysteresis loops as $\omega$ goes to zero in (i) to (iii).
\end{exercise} 

\begin{exercise}\label{exercisenonlinear}  For $y, u\in \mathbb R$, consider
\begin{subequations}\label{eqSONODE35}
\begin{align}
\ddot y(t) &=-\dot y(t) - 20y(t)^3((y(t))-0.3)(y(t)+0.5) +u(t),\\
y(0) &= -0.56, \dot y(0) =0.
\end{align}
\end{subequations} 
(a) Rewrite \eqref{eqSONODE35} into a system of first-order equations.\\
(b) Find the equilibria of \eqref{eqSONODE35}, and then determine their stability. \\
(c) For constant values of $u(t),$ find the values of $u$ for which \eqref{eqSONODE35} changes stability, and plot the corresponding stability diagram. \\
(d) For $u(t)=\sin(\omega t)$, modify the MATLAB function in exercise~\ref{exercisecode} to demonstrate the presence of hysteresis for various values of $\omega$. 
\end{exercise}

\begin{exercise} Repeat (a) to (d) in Exercise~\ref{exercisenonlinear} for 
\begin{subequations}\label{eqSONODE15}
\begin{align}
\ddot y(t) &=- 5\dot y(t) - (y(t))^2(y(t)+1)+u(t),\\
y(0) &= -1, \dot y(0) =0.
\end{align}
\end{subequations}
\end{exercise}

\bibliographystyle{siamplain}
\bibliography{references}

\end{document}